\documentclass[moor,sglanonrev]{informs4}
\RequirePackage{endnotes}

\SingleSpacedXII

\usepackage{tikz}
\usepackage{algorithmic}
\usepackage{array}
\usepackage{amsmath,epsfig,comment}
\usepackage{tcolorbox}
\usepackage{amssymb}
\usepackage{graphicx}
\usepackage{amstext}
\usepackage{mathrsfs}
\usepackage{multirow}
\usepackage{threeparttable}
\usepackage{subfigure}
\usepackage[multiple]{footmisc}
\usepackage{color}

\usepackage{amssymb,amsmath,cite}
\usepackage{epsfig}
\usepackage{color}
\usepackage{bm}
\usepackage[ruled,vlined,linesnumbered]{algorithm2e}
\usepackage{graphicx}
\usepackage{subfigure}
\usepackage{algorithmic}
\usepackage{multirow}
\usepackage{cite}
\usepackage{hyperref}
\usepackage{tikz}
\usetikzlibrary{arrows.meta,positioning}
\hypersetup{
	colorlinks=true,
	linkcolor=black,
	filecolor=magenta,
	urlcolor=cyan,
	citecolor=blue
}
\usepackage{extarrows}

\makeatletter
\def\@thmcountersep{.}
\makeatother

\usepackage{float}
\usepackage{booktabs}
\usepackage{multirow}

\DeclareMathAlphabet\mathbfcal{OMS}{cmsy}{b}{n}

\newcommand{\rgrad}{\mathrm{grad}\,}
\newcommand{\prox}{\mathrm{prox}}
\newcommand{\proj}{\mathrm{proj}}
\usepackage[sort&compress]{natbib}
 \bibpunct[, ]{[}{]}{,}{n}{}{,}%
 \def\bibfont{\small}%
\EquationsNumberedThrough    

\TheoremsNumberedThrough     
\ECRepeatTheorems  %

\MANUSCRIPTNO{MOOR-0001-2024.00}

\begin{document}


\RUNAUTHOR{Xu, Jiang, Liu, and So}

\RUNTITLE{A RADA Algorithmic Framework for Nonconvex-Linear Minimax Problems on Riemannian Manifolds}

\TITLE{A Riemannian Alternating Descent Ascent Algorithmic Framework for Nonconvex-Linear Minimax Problems on Riemannian Manifolds}

\ARTICLEAUTHORS{%
\AUTHOR{Meng Xu}
\AFF{ICMSEC, Academy of Mathematics and Systems Science, Chinese Academy of Sciences, and University of Chinese Academy of Sciences, \EMAIL{xumeng22@mails.ucas.ac.cn}}

\AUTHOR{Bo Jiang}
\AFF{Corresponding author. Ministry of Education Key Laboratory of NSLSCS, School of Mathematical Sciences, Nanjing Normal University, \EMAIL{jiangbo@njnu.edu.cn}}

\AUTHOR{Ya-Feng Liu}
\AFF{Corresponding author. Ministry of Education Key Laboratory of Mathematics
	and Information Networks, School of Mathematical Sciences, Beijing
	University of Posts and Telecommunications,  \EMAIL{yafengliu@bupt.edu.cn}}

\AUTHOR{Anthony Man-Cho So}
\AFF{Department of Systems Engineering and Engineering Management, The Chinese University of Hong Kong, \EMAIL{manchoso@se.cuhk.edu.hk}}
} 

\ABSTRACT{%
In this paper, we consider a class of nonconvex-linear minimax problems on Riemannian manifolds, which find wide applications in machine learning and signal processing. For solving this class of problems, we develop a flexible Riemannian alternating descent ascent (RADA) algorithmic framework. Within this framework, we propose two easy-to-implement yet efficient algorithms that alternately perform one or multiple projected/Riemannian gradient descent steps and a proximal gradient ascent step at each iteration.  We show that the proposed RADA algorithmic framework can find both an $\varepsilon$-Riemannian-game-stationary point and  an $\varepsilon$-Riemannian-optimization-stationary point within $\mathcal{O}(\varepsilon^{-3})$ iterations, achieving the best-known iteration complexity. We also reveal intriguing similarities and differences between the algorithms developed within our proposed framework and existing algorithms, thus providing important insights into the improved efficiency of the former. Lastly, we present numerical results on sparse principal component analysis (PCA), fair PCA, and sparse spectral clustering to demonstrate the superior performance of the proposed algorithms.
}%

\FUNDING{Meng Xu and Ya-Feng Liu are supported in part by the National Natural Science Foundation of China (NSFC) under Grant 12371314 and Grant 12021001. 
Bo Jiang is supported by the National Natural Science Foundation of China (NSFC) under Grant 12522116 and Grant 12371314.
Anthony Man-Cho So is supported in part by the Hong Kong Research Grants Council (RGC) General Research Fund (GRF) project CUHK 14204823.
}




\maketitle


\section{Introduction}\label{sec:introduction}
In this paper, we consider Riemannian nonconvex-linear (NC-L) minimax problems of the form
\begin{equation}\label{prob:p1}
	\min_{x\in\mathcal{M}}  \, \max_{y\in\mathcal{E}_2}\,\left\{F(x,y):=f(x)+\langle \mathcal{A}(x),\,y\rangle-h(y)\right\},
\end{equation}
where $\mathcal{M}$ is a Riemannian manifold embedded in a finite-dimensional Euclidean space $\mathcal{E}_1$, $\mathcal{E}_2$ is another finite-dimensional Euclidean space,  $f:\mathcal{E}_1\to\mathbb{R}$ is a continuously differentiable function, $\mathcal{A}:\mathcal{E}_1\to\mathcal{E}_2$ is a smooth mapping, and $h:\mathcal{E}_2\to(-\infty,+\infty]$ is a proper closed convex function with a compact domain and a tractable proximal mapping.  Many applications in machine learning and signal processing give rise to instances of problem \eqref{prob:p1}. We briefly introduce three such applications below and refer the interested reader to \cite{chen2020proximal,hu2020brief,li2021weakly,xiao2021exact,huang2022riemannian,li2023riemannian,liu2024survey} for more applications.  
\subsection{Motivating Applications}\label{subsec: application}
In the following, we use $A=[a_1,\,a_2,\,\ldots,\,a_N]\in\mathbb{R}^{d\times N}$ to denote a data matrix, where each of the $N$ columns corresponds to a data sample with $d$ attributes.

\begin{itemize}	
	\item[1.] \textbf{Sparse principal component analysis (SPCA).}
	Each principal component of $A$ obtained through the classic PCA is a linear combination of $d$ attributes, making it difficult to interpret the derived principal components as new features \cite{zou2018selective}. To achieve a good balance between dimension reduction and interpretability, SPCA seeks principal components with very few nonzero components.  This motivates the following formulation of SPCA \cite{jolliffe2003modified}:
	\begin{equation}\label{prob: SPCA}
		\min_{X\in\mathcal{S}(d,r)} \left\{-\langle AA^\top,\,XX^\top\rangle+\mu\|X\|_1\right\}.
	\end{equation}
	Here,  $\mathcal{S}(d,r) = \{X \in \mathbb{R}^{d \times r}\mid X^\top X = I_r\} $ is the Stiefel manifold with $I_r$ being the $r$-by-$r$ identity matrix, $\mu>0$ is the weighting parameter, and $\|X\|_1=\sum_{i,j}|X_{ij}|$ is the $\ell_1$-norm of the matrix $X$.  By Fenchel duality, the SPCA problem \eqref{prob: SPCA} can be equivalently reformulated as
	\begin{equation}\label{prob: SPCA ref minmax}
		\min_{X\in\mathcal{S}(d,r)}\max_{Y\in\mathcal{Y}} \left\{-\langle AA^\top,\,XX^\top\rangle + \langle Y,\,X\rangle\right\},
	\end{equation}
	where $\mathcal{Y}=\{Y\in\mathbb{R}^{d\times r}\mid\|Y\|_\infty\leq \mu \}$ with $\|Y\|_\infty=\max_{i,j}\,|Y_{ij}|$. Problem \eqref{prob: SPCA ref minmax} is an instance of problem \eqref{prob:p1} with $\mathcal{M}$ being the Stiefel manifold, $\mathcal{A}(X)=X$, and $h$ being the indicator function of the compact set $\mathcal{Y}$.
	
	\item[2.] \textbf{Fair principal component analysis (FPCA).} Suppose that the $N$ data samples belong to $ m $ groups according to certain clustering, and each group $i$ corresponds to a $d \times n_i$ submatrix $A_i$ with $\sum_{i = 1}^m n_i = N$. Recall that classic PCA aims to find a subspace with dimension $r < d$ to minimize the total reconstruction error, or equivalently, to maximize the total variance. As such, certain groups may suffer a higher reconstruction error than the others \cite{samadi2018price}.  To reduce such disparity among different groups, FPCA minimizes the maximum reconstruction error among the $m$ groups, which is equivalent to maximizing the minimum variance among the $m$ groups. Mathematically, the FPCA problem can be formulated as \cite{samadi2018price,zalcberg2021fair,xu2023efficient2,shen2025hidden}  
	\begin{equation*}\label{prob:FPCA}
		\min_{X\in\mathcal{S}(d,r)}\max_{i\in\{1,\,2,\,\ldots,\,m\}}\ -\langle A_iA_i^\top,\,XX^\top\rangle,
	\end{equation*}
	which is equivalent to
	\begin{equation}\label{prob: FPCA ref minmax}
		\min_{X\in\mathcal{S}(d,r)}\max_{y \in \Delta_m}\, -\sum_{i=1}^{m}y_i\langle A_iA_i^\top,\,XX^\top\rangle.
	\end{equation}
	Here, $ \Delta_m=\{y\in\mathbb{R}^m\,|\,\sum_{i=1}^{m}y_i=1,\,y_i\geq0,\,i=1,\,2,\,\ldots,\,m\}$ is the standard simplex in $\mathbb{R}^m$. Note that problem \eqref{prob: FPCA ref minmax} is an instance of problem \eqref{prob:p1}, where $\mathcal{M}$ is the Stiefel manifold, $\mathcal{A}: \mathbb{R}^{d \times r} \rightarrow \mathbb{R}^m$ is given by $\mathcal{A}(X)_i = -\langle A_i A_i^\top,\,XX^\top \rangle$ for $i=1,2,\ldots,m$, and $h$ is the indicator function of $\Delta_m$.
	
	\item[3.] \textbf{Sparse spectral clustering (SSC).} 
	This task aims to divide $N$ data samples into $m$ groups, each of which consists of similar data points.  Spectral clustering constructs a symmetric affinity matrix $W= [W_{ij}]_{N \times N}$, where $W_{ij}\geq 0$ measures the pairwise similarity between two samples $a_i$ and $a_j$.  To promote the sparsity and interpretability of spectral clustering, the works \cite{lu2016convex,lu2018nonconvex,wang2022manifold} proposed SSC, which entails solving the following problem:
	\begin{equation}\label{probSSC}
		\min_{Q\in\mathcal{G}(N,m)} \left\{\langle L,\,Q\rangle + \mu\|Q\|_1\right\}.
	\end{equation}
	Here, $\mathcal{G}(N, m)= \{XX^\top  \in \mathbb{R}^{N \times N} \mid X\in \mathcal{S}(N, m)\}$ is the Grassmann manifold \cite{bendokat2024grassmann}, which can be regarded as a submanifold in the Euclidean space \cite{sato2014optimization}; $L= I_N - S^{-1/2} W S^{-1/2}$ is the normalized Laplacian matrix with $S^{1/2}$ being the diagonal matrix with diagonal elements $\sqrt{s_1},\,\sqrt{s_2},\,\ldots,\,\sqrt{s_N}$ and $s_i = \sum_j W_{ij}$; and $\mu>0$ is the weighting parameter. By Fenchel duality,  problem \eqref{probSSC} can be equivalently reformulated as 
	\begin{equation}\label{prob: SSC ref minimax}
		\min_{Q\in\mathcal{G}(N,m)} \max_{Y\in\mathcal{Y}} \left\{\langle L,\,Q\rangle + \langle Y,\,Q\rangle\right\}.
	\end{equation}
	Again, problem \eqref{prob: SSC ref minimax} is an instance of problem \eqref{prob:p1} with $\mathcal{M}$ being the Grassmann manifold, $\mathcal{A}(Q)=Q$, and $h$ being the indicator function of the compact set {$\mathcal{Y}:=\{Y \in \mathbb{R}^{N \times N}\mid \|Y\|_{\infty} \leq \mu\}$}.
\end{itemize}

\subsection{Related Works}
Minimax problems in the Euclidean space---i.e., the setting where $\mathcal{M} \subseteq \mathcal{E}_1$ is an arbitrary closed convex set in problem \eqref{prob:p1}---have been the focus of many recent studies. For general smooth nonconvex-concave (NC-C) minimax problems, many nested-loop algorithms have been developed and analyzed \cite{nouiehed2019solving,thekumparampil2019efficient,lin2020near,rafique2022weakly,yang2020catalyst,kong2021accelerated,ostrovskii2021efficient}. 
Nevertheless, single-loop algorithms have gained growing interest due to their simplicity. A classic single-loop algorithm is the gradient descent ascent (GDA) method, which performs a (projected) gradient descent step on $x$ and a (projected) gradient ascent step on $y$ at each iteration. For smooth nonconvex-strongly concave (NC-SC) minimax problems, GDA can find an $\varepsilon$-stationary point within $\mathcal{O}(\varepsilon^{-2})$ iterations \cite{lin2020gradient}. However, when $F(x,\cdot)$ is not strongly concave for some $x$, GDA may encounter oscillations even for bilinear games. 
Many improved GDA variants have been proposed to overcome this issue \cite{jin2020local,lu2020hybrid,zhang2020single,yang2022faster,xu2023unified,li2022nonsmooth,wu2023efficient,xu2024derivative,lin2020gradient}.  One powerful approach is to design an appropriate regularized objective function. Such an approach was used by Xu et al. \cite{xu2023unified} and Zhang et al. \cite{zhang2020single} to develop the alternating gradient projection algorithm and the smoothed GDA algorithm, respectively. Both algorithms can find an $\varepsilon$-stationary point of a general smooth NC-C minimax problem within $\mathcal{O}(\varepsilon^{-4})$ iterations.  Later, Li et al. \cite{li2022nonsmooth} extended the techniques in \cite{zhang2020single} and developed the smoothed proximal linear descent ascent algorithm. The algorithm enjoys, among other things, an iteration complexity of $\mathcal{O}(\varepsilon^{-4})$ for finding  an $\varepsilon$-stationary point of certain nonsmooth NC-C minimax problem.  For the special case of NC-L minimax problems, Pan et al. \cite{pan2021efficient}, Shen et al. \cite{shen2023zeroth}, and He et al. \cite{he2024approximation} used the regularization approach to develop algorithms that can reach $\varepsilon$-stationarity within $\mathcal{O}(\varepsilon^{-3})$ iterations, which is the best complexity known to date.

Recently, there has been growing interest in minimax problems on Riemannian manifolds, but existing works on the topic are few. These works can be divided into two groups based on the nature of the constraints: (i) The constraints on $x$ and $y$ are manifolds \cite{zhang2023sion,jordan2022first,han2023riemannian,han2023nonconvex} and (ii) the constraint on $x$ is a manifold while the constraint on $y$ is a closed convex set \cite{huang2023gradient,xu2023efficient2}.  Among the works in the first group, Zhang et al. \cite{zhang2023sion} generalized the classic Sion's minimax theorem to {geodesic metric spaces} and proposed a Riemannian corrected extragradient (RCEG) algorithm for solving geodesically convex-geodesically concave (GC-GC) minimax problems. Jordan et al. \cite{jordan2022first} analyzed both the last-iterate and average-iterate convergence of RCEG and a Riemannian gradient descent ascent (RGDA) algorithm in various GC-GC settings.  For problems without geodesic convexity, Han et al. \cite{han2023riemannian} extended the Hamiltonian gradient method to the Riemannian setting and established its global linear convergence under the assumption that the Riemannian Hamiltonian of $F$ satisfies the Riemannian Polyak-\L ojasiewicz (P\L) condition.  In addition, Han et al. \cite{han2023nonconvex} proposed several second-order methods that are proven to asymptotically converge to local minimax points of nonconvex-nonconcave minimax problems under strong assumptions.  As for the works in the second group,  Huang and Gao \cite{huang2023gradient} proposed a RGDA algorithm different from the one in \cite{jordan2022first} and established its iteration complexity of $\mathcal{O}(\varepsilon^{-2})$ for finding an $\varepsilon$-stationary point of an NC-SC minimax problem. 
Later, the work \cite{xu2023efficient2} presented an alternating Riemannian/projected gradient descent ascent (ARPGDA) algorithm with an iteration complexity of $\mathcal{O}(\varepsilon^{-3})$ for finding an $\varepsilon$-stationary point of the Riemannian NC-L minimax problem \eqref{prob:p1} when $h$ is the indicator function of a convex compact set. 

From the above discussion, we observe that existing works on minimax problems over Riemannian manifolds do not address the settings of problem \eqref{prob:p1} that we are interested in. Indeed, neither the geodesic convexity assumption in \cite{zhang2023sion,jordan2022first} nor the strong concavity assumption in \cite{huang2023gradient} apply to the motivating applications introduced in Section \ref{subsec: application}. Additionally, the Riemannian P\L  ~condition of the Hamiltonian function in \cite{han2023riemannian} is difficult to verify and is not known to hold for these applications. On the other hand, the second-order methods in \cite{han2023nonconvex} cannot be used to tackle problem \eqref{prob:p1} since the Hessian of $F(x,\cdot)$ is singular.  Given the above background, \emph{we are motivated to develop simple yet efficient single-loop first-order methods with strong theoretical guarantees for solving the Riemannian NC-L minimax problem \eqref{prob:p1}}. 

One possible approach is based on the following reformulation of \eqref{prob:p1} as a Riemannian nonsmooth composite problem:
\begin{equation}\label{prob:Riemannian nonsmooth composite}
	\min_{x\in\mathcal{M}}   \left\{f(x)+h^*(\mathcal{A}(x))\right\}.
\end{equation}
Here, $ h^*$ is the conjugate function of $h$.  Indeed, since $h$ is assumed to have a compact domain, we know from \cite[Theorem 4.23]{beck2017first} that $h^*$ is Lipschitz continuous.  Then, we can apply a host of algorithms to tackle problem \eqref{prob:Riemannian nonsmooth composite}, including  Riemannian subgradient-type methods \cite{borckmans2014riemannian,hosseini2017riemannian,hosseini2018line,li2021weakly,hu2023constraint}, Riemannian proximal gradient-type methods \cite{chen2020proximal,huang2022riemannian,huang2023inexact,wang2022manifold,chen2024nonsmooth,liu2024penalty}, Riemannian smoothing-type algorithms \cite{beck2023dynamic,peng2023riemannian,zhang2023riemannian}, and splitting-type methods \cite{lai2014splitting,kovnatsky2016madmm,deng2023manifold,li2023riemannian,zhou2023semismooth,deng2024oracle,liu2020simple}. 
In particular, when $\mathcal{A}$ is the identity map $\mathcal{I}$ and $h^*$ is a real-valued weakly convex function, 
Li et al. \cite{li2021weakly} proposed a family of Riemannian subgradient (RSG) methods and showed that all of them have an iteration complexity of $\mathcal{O}(\varepsilon^{-4})$ for getting certain stationarity measure of problem \eqref{prob:Riemannian nonsmooth composite} below $\varepsilon$. When $\mathcal{A}=\mathcal{I}$, Chen et al. \cite{chen2020proximal} proposed the manifold proximal gradient (ManPG) method with an outer iteration complexity of $\mathcal{O}(\varepsilon^{-2})$ for finding an $\varepsilon$-stationary point of problem \eqref{prob:Riemannian nonsmooth composite}.  Here, the outer iteration complexity refers to the number of subproblems that need to be solved in the nested-loop algorithms. In particular, the complexity of solving each subproblem is not counted in the outer iteration complexity.
When $\mathcal{A}$ is a linear mapping, Beck and Rosset \cite{beck2023dynamic} proposed a dynamic smoothing gradient descent method (DSGM) with an iteration complexity of $\mathcal{O}(\varepsilon^{-3})$ for finding an $\varepsilon$-stationary point of problem \eqref{prob:Riemannian nonsmooth composite}.  Using the splitting technique, Li et al. \cite{li2023riemannian} proposed a Riemannian ADMM (RADMM) with an iteration complexity of $\mathcal{O}(\varepsilon^{-4})$ and Deng et al. \cite{deng2024oracle} proposed a manifold inexact augmented Lagrangian (ManIAL) method with an iteration complexity of $\mathcal{O}(\varepsilon^{-3})$ for finding an $\varepsilon$-stationary point of problem \eqref{prob:Riemannian nonsmooth composite}. When $\mathcal{A}$ is a nonlinear mapping, Wang et al. \cite{wang2022manifold} proposed a manifold proximal linear (ManPL) method and established its outer iteration complexity of $\mathcal{O}(\varepsilon^{-2})$ for finding an $\varepsilon$-stationary point of problem \eqref{prob:Riemannian nonsmooth composite}.

In this paper, we propose an alternative approach, which is based on exploiting the minimax structure of problem \eqref{prob:p1} and can also be used to solve problem \eqref{prob:Riemannian nonsmooth composite}. As seen from Table \ref{tab: complexity comparision}, which summarizes the applicability and complexity of different methods, our approach is able to address more general problem settings and has the best complexity known to date among single-loop algorithms for solving problem \eqref{prob:p1}. Moreover, as we shall see in Section \ref{sec: connections}, our approach reveals intriguing similarities and differences between our proposed algorithms for problem \eqref{prob:p1} and existing algorithms for problem \eqref{prob:Riemannian nonsmooth composite}.

\begin{table}[t]
	\fontsize{10.5pt}{\baselineskip}\selectfont
	\centering
	\renewcommand{\arraystretch}{1}
	\tabcolsep=0.05cm
	\caption{A summary of methods with iteration complexity results for solving problem \eqref{prob:p1} or \eqref{prob:Riemannian nonsmooth composite}, where “$\star$” denotes the outer iteration complexity. The results for HiBSA (developed for general Euclidean NC-C minimax problems) and AGP (developed for Euclidean NC-L minimax problems) are stated for the Euclidean NC-L minimax setting.
	}
	\begin{tabular}{c|c|c|c|c}
		\hline
		{Algorithm} & {Mapping $\mathcal{A}$} & {Function $h^*$} & {Complexity}& Type\\ 
		\hline
		RSG \cite{li2021weakly} & $\mathcal{I}$& real-valued weakly convex&  $\mathcal{O}(\varepsilon^{-4})$ & single-loop \\
		\hline
		{ManPG} \cite{chen2020proximal} & $\mathcal{I}$ & convex, Lipschitz continuous & $\mathcal{O}(\varepsilon^{-2})^\star$&nested-loop\\
		\hline
		ManIAL \cite{deng2024oracle} & linear & convex, Lipschitz continuous &$\mathcal{O}(\varepsilon^{-3})$&nested-loop \\
		\hline
		{DSGM} \cite{beck2023dynamic} & linear & convex, Lipschitz continuous& $\mathcal{O}(\varepsilon^{-3})$& single-loop\\
		\hline
		RADMM \cite{li2023riemannian} & linear & convex, Lipschitz continuous & $\mathcal{O}(\varepsilon^{-4})$ & single-loop \\
		\hline
		ManPL \cite{wang2022manifold} & nonlinear & convex, Lipschitz continuous& $\mathcal{O}(\varepsilon^{-2})^\star$& nested-loop\\
		\hline
		ARPGDA \cite{xu2023efficient2} & nonlinear & support function of a convex compact set& $\mathcal{O}(\varepsilon^{-3})$ & single-loop \\
		\hline
		HiBSA \cite{lu2020hybrid} & nonlinear & support function of a convex compact set& $\mathcal{O}(\varepsilon^{-4})$ & single-loop \\
		\hline
		AGP \cite{pan2021efficient} & nonlinear & support function of a convex compact set& $\mathcal{O}(\varepsilon^{-3})$ & single-loop \\
		\hline
		RADA-PGD/RGD [this paper] & nonlinear  & convex, Lipschitz continuous & $\mathcal{O}(\varepsilon^{-3})$ & single-loop \\
		\hline \end{tabular}
	\label{tab: complexity comparision}
\end{table}
\subsection{Our Contributions}
We now summarize the main contributions of this paper as follows.
\begin{itemize}
	\item[$\bullet$] We propose a flexible Riemannian alternating descent ascent (RADA) algorithmic framework for solving problem \eqref{prob:p1}. At the heart of our framework is the observation that the value function associated with the maximization of certain regularized version of the objective function $F$ with respect to $y$ is smooth and has an easily computable gradient; see \eqref{potential Phi_k} and \eqref{nabla Phik}. At each iteration, the proposed RADA algorithmic framework finds an approximate minimizer of the said value function to update $x$ and performs one proximal gradient ascent step to update $y$.  The value function contains information about the best response of $y$, namely, the maximization of the objective function $F(x,y)$ with respect to $y$ for any given $x$. Therefore, minimizing the proposed value function can be regarded as approximately solving a surrogate of problem \eqref{prob:p1}. It should be noted that our approach to updating $x$ is different from simply finding an (approximate) minimizer of some surrogate of $F(\cdot,y)$ for a fixed $y$, as is done in, e.g., \cite{xu2023unified,pan2021efficient,lu2020hybrid,he2024approximation}, and is more advantageous in practice. Within the RADA algorithmic framework, we provide two customized simple yet efficient single-loop first-order algorithms, namely, RADA-PGD and RADA-RGD, which perform a fixed number of projected gradient descent steps and Riemannian gradient descent steps to update $x$ at each iteration, respectively. Interestingly, even when adapted to the setting where $\mathcal{M}$ is simply a closed convex set in a Euclidean space, our proposed RADA algorithmic framework, along with the two customized single-loop algorithms, is different from existing methods for solving Euclidean NC-L minimax problems (such as those in \cite{xu2023unified,pan2021efficient,lu2020hybrid,he2024approximation}) and is, to the best of our knowledge, new.
	\item[$\bullet$] We prove that our proposed RADA algorithmic framework can find an $\varepsilon$-Riemannian-game-stationary ($\varepsilon$-RGS) point and an $\varepsilon$-Riemannian-optimization-stationary ($\varepsilon$-ROS) point of problem \eqref{prob:p1} within $\mathcal{O}(\varepsilon^{-3})$ iterations. This iteration complexity result matches the best-known complexity result for general NC-L minimax problems where the variable $x$ is constrained on a nonempty compact convex set instead of a Riemannian manifold \cite{pan2021efficient, he2024approximation}. It also matches the best-known complexity result for the Riemannian nonsmooth composite problem \eqref{prob:Riemannian nonsmooth composite} with a linear mapping $\mathcal{A}$ \cite{beck2023dynamic,deng2024oracle}.  Furthermore, while one typically needs extra computation to convert between an $\varepsilon$-game stationary point and an $\varepsilon$-optimization stationary point in the Euclidean setting (see, e.g., \cite{lin2020gradient,yang2022faster}), it is worth noting that using our RADA algorithmic framework, the complexity of finding an $\varepsilon$-RGS point is of the same order as that of finding an $\varepsilon$-ROS point.
	\item[$\bullet$] We elaborate on the relationship between (i) the RADA algorithmic framework and RADA-RGD in particular, which are designed for problem \eqref{prob:p1}, and (ii) the Riemannian augmented Lagrangian method (RALM) \cite{deng2023manifold,zhou2023semismooth,deng2024oracle} and RADMM \cite{li2023riemannian}, which are designed for problem \eqref{prob:Riemannian nonsmooth composite}. Roughly speaking, the subproblem for updating $x$ in RALM is the same as that in RADA when certain regularization parameter in the aforementioned value function is zero.  However, the stopping criterion for the subproblem in RADA is less stringent than that in RALM. Moreover, we show that RADA-RGD with certain parameter setting gives rise to a new Riemannian symmetric Gauss-Seidel ADMM (sGS-ADMM). These similarities and differences provide important insights into understanding the superior performance of our proposed algorithms over existing ones \cite{zhou2023semismooth,li2023riemannian}.  
	\item[$\bullet$] We present numerical results on SPCA, FPCA, and SSC to demonstrate the advantages of our proposed algorithms over existing state-of-the-art algorithms for the corresponding problems. 	Our numerical results on SPCA show that when compared against the nested-loop algorithms ManPG \cite{chen2020proximal} and RALM \cite{zhou2023semismooth}, which employ second-order methods to solve their respective subproblems, our proposed RADA-RGD can return a solution of comparable quality with much lower computational cost. In addition, our numerical results on SSC and FPCA show that our proposed algorithms outperform several state-of-the-art single-loop first-order algorithms \cite{beck2023dynamic,xu2023efficient2,lu2018nonconvex,li2023riemannian} in terms of both solution quality and speed.
\end{itemize}

\subsection{Organization}
The rest of the paper is organized as follows.  We first introduce the notation and cover some preliminaries on Riemannian optimization and convex analysis in Section \ref{sec:preliminaries}.  We then present our proposed algorithmic framework and establish its iteration complexity in Section \ref{sec:proposed algorithm framework}.  Subsequently, we develop two customized algorithms within the proposed framework in Section \ref{sec: two speicific}.  We reveal the connections between the proposed algorithms and some existing ones in Section \ref{sec: connections}. Then, we present numerical results on SPCA, FPCA, and SSC to illustrate the efficiency of the proposed algorithms in Section \ref{sec: numerical results}. Finally, we draw some conclusions in Section \ref{sec: concluding remarks}.

\section{Notation and Preliminaries}\label{sec:preliminaries}
We first introduce the notation and some basic concepts in Riemannian optimization \cite{absil2008optimization,boumal2023introduction}. 
We use $\langle \,\!\cdot\,, \cdot\rangle$ and $\| \cdot \|$ to denote the standard inner product and its induced norm on the Euclidean space $\mathcal{E}$, respectively. 
Let $\mathcal{X}$ be a subset of $\mathcal{E}$. We use $\mathrm{dist}(y,\,\mathcal{X})=\inf_{x\in\mathcal{X}}\,\|y-x\|$ to denote the distance from $y\in\mathcal{E}$ to  $\mathcal{X}$ and $\mathrm{conv}\,\mathcal{X}$ to denote the convex hull of $\mathcal{X}$.
For a linear mapping $\mathcal{A}:\mathcal{E}_1\to\mathcal{E}_2$, we use $\mathcal{A}^\top:\mathcal{E}_2\to\mathcal{E}_1$ to denote its adjoint mapping.
Let $\mathcal{M}$ be a submanifold embedded in $\mathcal{E}$ and $ \mathrm{T}_{x}\mathcal{M} $ denote the tangent space to $\mathcal{M}$ at $x\in\mathcal{M}$.
Throughout this paper, we take the standard inner product $\langle \cdot,\cdot \rangle$ on the Euclidean space $\mathcal{E}$ as the Riemannian metric on $\mathcal{M}$. Then, for a smooth function $f: \mathcal{E} \rightarrow \mathbb{R}$ and a point $x \in \mathcal{M}$, the Riemannian gradient of $f$ at $x$ is given by
$\rgrad f(x)=\proj_{\mathrm{T}_{x}\mathcal{M}}(\nabla f(x))$, where $\proj_{\mathcal{X}}(\cdot)$ is the Euclidean projection operator onto a closed set $\mathcal{X}$ and $\nabla f(x)$ is the Euclidean gradient of $f$ at $x$. 
A retraction at $x \in \mathcal{M}$ is a smooth mapping $\mathrm{R}_x: \mathrm{T}_x \mathcal{M} \to \mathcal{M}$ satisfying (i) $\mathrm{R}_x(\mathbf{0}_x) = x$, where $\mathbf{0}_x$ is the zero element in $\mathrm{T}_x \mathcal{M}$; (ii) $\frac{\mathrm{d}}{\mathrm{d} t} \mathrm{R}_x (t v)|_{t = 0} = v$ for all $v \in \mathrm{T}_x \mathcal{M}$. 

Now, we review some basic notions in convex analysis \cite{rockafellar2009variational,beck2017first}. 
Let $g:\mathcal{E}\to(-\infty,+\infty]$ be a proper closed convex function whose domain is given by $\mathrm{dom}\ g:= \{x \in \mathcal{E} \mid g(x) < +\infty\}$. Its conjugate function is defined as
$
	g^*(y):=\max_{x\in\mathcal{E}}\left\{\langle y,\,x\rangle-g(x)\right\}.
$
The function $g$ is called $\lambda$-strongly convex if $g(\cdot)-\frac{\lambda}{2}\|\cdot\|^2$ is convex.
For a given constant $\lambda > 0$, the proximal mapping and the Moreau envelope of $g$ are defined as 
\begin{equation*}
	\prox_{ \lambda g}(x)=\argmin_{u\in\mathcal{E}}\left\{g(u)+\frac{1}{2\lambda}\|x-u\|^2\right\}
\end{equation*} 
and
\begin{equation*}\label{Moreau envelope}
	M_{{\lambda g}}(x)=\min_{u\in\mathcal{E}}\left\{g(u)+\frac{1}{2\lambda}\|x-u\|^2\right\},
\end{equation*}
respectively. 
For any $x \in \mathcal{E}$, we have the following so-called Moreau decomposition and Moreau envelope decomposition:
\begin{align}
	&x=\prox_{\lambda g}(x)+\lambda\prox_{{g^*/\lambda}}\left(\frac{x}{\lambda}\right),\label{Moreau Decomposition}\\
	&\frac{1}{2\lambda}\|x\|^2=M_{{\lambda g}}(x)+M_{{g^*/\lambda}}\left(\frac{x}{\lambda}\right).\label{Moreau envelope decomposition}
\end{align}
The following theorem characterizes the smoothness of the Moreau envelope.
\begin{theorem}\label{gradient of Moreau}\textbf{(\hspace{-0.01cm}\cite[Theorems 6.39 and 6.60]{beck2017first})}
	Let $g:\mathcal{E}\to(-\infty,+\infty]$ be a proper closed convex function and $\lambda>0$. Then, for any $x\in\mathcal{E}$,
	\begin{equation}\label{gradient of Moreau envelope}
		\nabla M_{\lambda g}(x)=\frac{1}{\lambda}(x-\prox_{\lambda g}(x))=\mathrm{prox}_{g^*/\lambda}\left(\frac{x}{\lambda}\right)\in\partial g(\prox_{\lambda g}(x)).
	\end{equation}
\end{theorem}
Moreover, $\nabla M_{\lambda g }$ is $\lambda^{-1}$-Lipschitz continuous.

\section{Proposed Algorithmic Framework}\label{sec:proposed algorithm framework}
In this section, we propose a flexible algorithmic framework for solving the Riemannian NC-L minimax problem \eqref{prob:p1} and establish its iteration complexity. 

\subsection{Proposed RADA Algorithmic Framework}\label{subsection: Algorithm framework}

A widely used idea for tackling the minimax problem \eqref{prob:p1} is to alternately update the variables $x$ and $y$ by solving an appropriate minimization problem for $x$ and an appropriate maximization problem for $y$. However, even a small perturbation in $x$ may cause a drastic change in the solution to the problem $\max_{y\in\mathcal{E}_2} F(x,y)$ \cite{sun2024dual}. Inspired in part by \cite{lu2020hybrid,pan2021efficient,he2024approximation,xu2023unified}, we address this issue by adding quadratic regularization terms to the objective function and defining the following value function at the $k$-th iteration:
\begin{equation}\label{potential Phi_k}
	\Phi_k(x):= \max_{y\in\mathcal{E}_2}\left\{F(x,y)-\frac{\lambda}{2}\|y\|^2-\frac{\beta_{k}}{2}\|y-y_{k}\|^2\right\}.
\end{equation}
Here, $\lambda > 0$ is the regularization parameter, $\beta_k \geq 0$ is the proximal parameter, and $y_k \in \mathrm{dom}\, h$ is the proximal center. Intuitively, the quadratic regularization term $\frac{\lambda}{2}\|y\|^2$ ensures the smoothness of $\Phi_k$ (as shown later in \eqref{nabla Phik}), while the proximal term $\frac{\beta_k}{2}\|y-y_k\|^2$ stabilizes the update of $y$. 

At the $k$-th iteration, our algorithmic framework computes the update $x_{k+1}$ by approximately solving the subproblem
\begin{equation}\label{subproblemphik}
	\min_{x\in\mathcal{M}}\,\Phi_k(x).
\end{equation}
Specifically, we need $x_{k+1}$ to satisfy the sufficient decrease condition
\begin{equation}\label{sufficient decrease}
	\Phi_k(x_{k+1})-\Phi_{k}(x_{k})\leq-\min\left\{c_\lambda\|\rgrad\Phi_{k}(x_k)\|,\,c'\right\}\|\rgrad\Phi_{k}(x_k)\|+\nu_k,
\end{equation}
where {$c_\lambda >0$ is a predetermined constant on the order of $\lambda$} (i.e., $c_{\lambda} = \mathcal{O}(\lambda)$), $c'>0$ is another predetermined constant, and 
\begin{equation*}\label{summable sequence}
	\{ \nu_k \} \in \mathbb{S}:=\left\{ \{\nu_{k}\}\;\middle|\;\sum_{k = 1}^{+\infty} \nu_k < +\infty,\,\nu_{k}\geq0,\,k=1,\,2,\,\ldots\right\}
\end{equation*}
is a nonnegative summable sequence. Then, we compute the update
\begin{equation}
	\begin{aligned}\label{AlgorithmupdateY}
		y_{k+1}& = \argmax_{y\in\mathcal{E}_2}\left\{F(x_{k+1},y)-\frac{\lambda}{2}\|y\|^2-\frac{\beta_{k}}{2}\|y-y_{k}\|^2\right\} \\ 
		& =\prox_{{h/(\lambda + \beta_k)}}\left(\frac{\mathcal{A}(x_{k+1}) + \beta_k y_k}{\lambda + \beta_k}\right).
	\end{aligned}	
\end{equation}
The above update of $y$ can also be regarded as performing a proximal gradient ascent step of the regularized function $F(x_{k+1},\cdot)-\frac{\lambda}{2}\|\cdot\|^2$. The proposed RADA algorithmic framework is formally presented in Algorithm \ref{Algorithmframework}.

\vspace{0.5cm}
\begin{algorithm}[t]
	\caption{RADA algorithmic framework for problem \eqref{prob:p1} }\label{Algorithmframework}
	Input $x_1\in\mathcal{M},\, y_1\in\mathrm{dom}\,h$, $\lambda>0$, $\beta_1\geq0$, $\rho>1$, {$c_\lambda>0$}, and $\{\nu_{k}\}\in\mathbb{S}$.
	
	\For{$k=1,\,2,\,\ldots$}
	{Approximately solve subproblem \eqref{subproblemphik} to find a point $x_{k+1}\in\mathcal{M}$ that satisfies the sufficient decrease condition \eqref{sufficient decrease}.
		
		Calculate $y_{k+1}$ via \eqref{AlgorithmupdateY}. 
		
		Update $\beta_{k+1}$ such that  $0\leq\beta_{k+1}\leq{\beta_{1}}/{(k+1)^\rho}$.
	}
\end{algorithm}
Several remarks on the proposed Algorithm \ref{Algorithmframework} are in order. First,
the requirement on the update $\beta_{k+1}$ is quite flexible. A simple approach is to set $\beta_{k}=\beta_{1}/k^\rho$ for all $k\ge1$. Alternatively, we can adaptively update $\beta_{k}$ for improved practical efficiency.
For instance, given constants $\tau_1\in (0,1)$ and $\tau_2 \in (0,1)$, we may set
\begin{equation}\label{RADA up beta}
	\beta_{k+1} = \frac{\beta_1^{(k+1)}}{(k+1)^\rho}\quad \mbox{with}\quad \beta_1^{(k+1)}=\left\{
	\begin{aligned}
		\tau_2\beta_{1}^{(k)}, &\quad \text {if } \delta_{k+1}\geq\tau_1\delta_k; \\ 
		\beta_{1}^{(k)},&\quad \text {if } \delta_{k+1}<\tau_1\delta_k,
	\end{aligned}
	\right.
\end{equation}
where $\beta_1^{(1)} = \beta_1$ and  $\delta_{k+1}$ is defined as
\begin{equation}\label{delta}
	{\delta_{k+1}} := \|(\lambda+\beta_k)y_{k+1}- \beta_{k}y_{k}\|_\infty, \quad\forall\,k\geq0.
\end{equation} 
Here, we adopt the convention that $\beta_0 = \beta_1$ and $y_0 = y_1$ when computing $\delta_1$.  
It is simple to see that the update $\beta_{k+1}$ in \eqref{RADA up beta} satisfies $0\leq\beta_{k+1}\leq{\beta_{1}}/{(k+1)^\rho}$ for all $k\geq1.$
As will be discussed at the end of Section \ref{Section relation RALM}, our proposed Algorithm \ref{Algorithmframework} is closely related to RALM. In particular, the parameter $\delta_{k+1}$ in \eqref{delta} can be viewed as a primal residual in RALM. 

Second, the sufficient decrease condition \eqref{sufficient decrease} is easy to satisfy, making the proposed algorithmic framework computationally efficient. For instance, we can find the desired $x_{k+1}\in\mathcal{M}$ by performing a simple Riemannian gradient descent step on $x$ with a suitable stepsize. 
To do this, we need to establish the smoothness of $\Phi_k$. By the fact that
\begin{align*} 
	&\langle \mathcal{A}(x),\,y\rangle-\frac{\lambda}{2}\|y\|^2-\frac{\beta_{k}}{2}\|y-y_{k}\|^2+\frac{\beta_{k}}{2}\|y_k\|^2\\
	=\,&\frac{1}{2(\lambda+\beta_{k})}\|\mathcal{A}(x)+\beta_{k}y_{k}\|^2 - \frac{\lambda+\beta_{k}}{2}\left\|y-\frac{\mathcal{A}(x)+\beta_{k}y_{k}}{\lambda+\beta_{k}}\right\|^2,
\end{align*}
we have from \eqref{potential Phi_k} and \eqref{Moreau envelope decomposition} that
\begin{equation}\label{Moreau Phi h} 
	\begin{aligned}
		\Phi_{k}(x)=&\, f(x)+\max_{y\in\mathcal{E}_2} \ \left\{\langle \mathcal{A}(x),\,y\rangle-h(y)-\frac{\lambda}{2}\|y\|^2-\frac{\beta_{k}}{2}\|y-y_{k}\|^2\right\}\\
		= &\,f(x) +\frac{1}{2(\lambda+\beta_{k})}\|\mathcal{A}(x)+\beta_{k}y_{k}\|^2- M_{{h/(\lambda+\beta_k)}}\left(\frac{\mathcal{A}(x)+\beta_{k}y_{k}}{\lambda+\beta_{k}}\right)-\frac{\beta_{k}}{2}\|y_k\|^2\\
		= &\,f(x) + M_{{(\lambda+\beta_k)h^*}}\left(\mathcal{A}(x)+\beta_{k}y_{k}\right)-\frac{\beta_{k}}{2}\|y_k\|^2.
	\end{aligned}
\end{equation}
Therefore, by Theorem \ref{gradient of Moreau}, the value function $\Phi_{k}$ is differentiable and \begin{equation}\label{nabla Phik}
	\nabla\Phi_{k}(x)=\nabla f(x)+\nabla\mathcal{A}(x)^\top\prox_{{h}/(\lambda+\beta_{k})}\left(\frac{\mathcal{A}(x) + \beta_k y_k}{\lambda + \beta_k}\right),
\end{equation} where $\nabla\mathcal{A}$ denotes the Jacobian mapping of $\mathcal{A}$.
More approaches for approximately solving the subproblem \eqref{subproblemphik} are presented in Section \ref{sec: two speicific}. The choice of the approach mainly depends on the structure of the underlying problem. The update $y_{k+1}$ in \eqref{AlgorithmupdateY} is also simple since $h$ is assumed to have a tractable proximal mapping. Therefore, our proposed algorithmic framework enjoys a low per-iteration computational cost. 

Finally, we highlight the difference between our proposed RADA algorithmic framework and some related algorithms in the literature.
The works \cite{pan2021efficient,lu2020hybrid,xu2023unified,he2024approximation} for solving Euclidean NC-L minimax problems minimize certain surrogate of $F(\cdot,y)$ with a fixed $y$ to update $x$, i.e.,
\begin{equation*}
	x_{k+1} \approx \argmin_{x \in \mathcal{M}} \left\{  { F(x,y_k)- \frac{\lambda}{2}\|y_k\|^2} \right\},
\end{equation*}
where $\mathcal{M}$ is a compact convex set.
By contrast, the proposed value function \eqref{potential Phi_k} at $x$ is the maximization of a regularized $F(x,\cdot)$, and the subproblem \eqref{subproblemphik} used to update $x_{k+1}$ is itself a minimax problem, i.e., 
\begin{equation*}
	 x_{k+1} \approx \argmin_{x \in \mathcal{M}}  \left\{{\max_{y\in\mathcal{E}_2}\left\{ F(x,y)- \frac{\lambda}{2}\|y\|^2 - \frac{\beta_{k}}{2}\|y-y_{k}\|^2 \right\}}\right\},
\end{equation*}
which can be viewed as a surrogate of problem \eqref{prob:p1}.
Thus, our proposed RADA is fundamentally different from the approaches in \cite{pan2021efficient,lu2020hybrid,xu2023unified,he2024approximation}.
To provide more insights into why the former could outperform the latter, we elaborate on the difference between RADA and the ARPGDA algorithm developed in \cite{xu2023efficient2}, which can be regarded as a Riemannian counterpart of the methods in \cite{pan2021efficient,lu2020hybrid,xu2023unified,he2024approximation}.
ARPGDA computes the update
\begin{align*}
	x_{k+1}=\mathrm{R}_{x_k}\left(-\zeta_k\,\mathrm{grad}_x\,F(x_k,y_k)\right),
\end{align*}
where $\mathrm{grad}_x\,F(x,y)$ denotes the Riemannian gradient of $F(x,y)$ with respect to $x$. Then, it computes the update $y_{k+1}$ according to \eqref{AlgorithmupdateY}.
Similar to the derivation of $\nabla\Phi_k(x)$, we have $\nabla_x F(x_k, y_k) = \nabla\Phi_{k-1}(x_k)$ and hence $$\mathrm{grad}_x\,F(x_k, y_k) = \rgrad \Phi_{k-1} (x_k).$$
Therefore, ARPGDA uses the Riemannian gradient of the value function $\Phi_{k-1}$ at $x_k$ to compute the update $x_{k+1}$. 
By contrast, our proposed framework uses the Riemannian gradient of the value function $\Phi_k$ at $x_k$ to compute the update $x_{k+1}$. In particular, our proposed framework makes use of the most recent information about $y$ (since $\Phi_k$ involves the latest iterate $y_k$) to compute the update $x_{k+1}$. This suggests that RADA could outperform ARPGDA. 

\subsection{Iteration Complexity}\label{sec: iteration complexity} 
In this subsection, we study the convergence behavior of the proposed RADA algorithmic framework. 
Towards that end, let us first introduce the stationarity measures that we are interested in. 
The following definition is motivated by the stationarity measures used in Euclidean minimax problems \cite{xu2023unified,li2022nonsmooth} and (Riemannian) nonsmooth composite problems \cite{beck2023dynamic,li2023riemannian,tian2024no}.
\begin{definition}\label{stationarypoint}
	Let $\varepsilon>0$ be some given constant.
	\begin{itemize}
		\item[(a)] The point $(x, y)\in\mathcal{M}\, \times \, \mathrm{dom}\, h$ is an $\varepsilon$-Riemannian-game-stationary ($\varepsilon$-RGS) point of problem (\ref{prob:p1}) if there exists a constant $\gamma>0$ such that
		\begin{equation*}\label{FNE2}
			\max\!\left\{\left\|\mathrm{grad}_x\,F(x, y)\right\|\hspace{-2pt},\,\frac{1}{\gamma}\left\|y-\prox_{\gamma h}(y+\gamma\mathcal{A}(x))\right\|\right\}\leq\varepsilon.
		\end{equation*}
		\item[(b)] The point $x\in\mathcal{M}$ is an $\varepsilon$-Riemannian-optimization-stationary ($\varepsilon$-ROS) point of problem \eqref{prob:p1} if there exists a point $p\in\mathcal{E}_2$ such that
		\begin{equation*}\label{Def optstationary}
			\max\!\left\{\mathrm{dist}\left(0,\,\rgrad f(x) + \proj_{\mathrm{T}_{x}\mathcal{M}}\left(\nabla \mathcal{A}(x)^\top \partial h^*(p)\right)\right)\hspace{-2pt},\,\left\|p-\mathcal{A}(x)\right\|\right\}\leq\varepsilon,
		\end{equation*}
		where $\partial h^*$ is the subdifferential of $h^*$.
	\end{itemize}
\end{definition}

For simplicity, denote
\begin{align}
	F_k(x,y)&:=F(x,y)-\frac{\lambda}{2}\|y\|^2-\frac{\beta_k}{2}\|y-y_{k}\|^2\label{F_k},\\
	y_{k+\frac{1}{2}}&:= \argmax_{y\in\mathcal{E}_2}\,F_k(x_k,y)=\prox_{{h}/(\lambda+\beta_{k})}\left(\frac{\mathcal{A}(x_k) + \beta_k y_k}{\lambda + \beta_k}\right),\label{y_k+1/2}\\
	\Phi(x)&:=\max_{y\in\mathcal{E}_2} F(x,y),\label{Phi}\\
	R&:=\max_{y\in \mathrm{dom}\, h}\,\|y\| <+\infty.\label{R}
\end{align}
Recalling that $\rho>1$, $0 \leq \beta_k \leq {\beta_1}/{k^\rho}$, and $\{ \nu_k \} \in \mathbb{S}$ in Algorithm \ref{Algorithmframework}, we have
\begin{equation}\label{Upsilon}
	\Upsilon:= \frac{\lambda R^2}{2}+\sum_{k=1}^{+\infty}(\nu_{k}+2\beta_{k}R^2)<+\infty.
\end{equation} 
We make the following blanket assumption.
\begin{assumption}\label{assumption: level bound}
	The function $\Phi$ defined in \eqref{Phi} is bounded from below on $\mathcal{M}$, namely, $\Phi^*:=\inf_{x \in \mathcal{M}} \Phi(x) > -\infty$.
\end{assumption}

Now, we are ready to establish the iteration complexity of the proposed RADA algorithmic framework for returning an $\varepsilon$-RGS point and an $\varepsilon$-ROS point of problem \eqref{prob:p1}. 
We start with the following estimate:
\begin{lemma}\label{lemma descent condition Phi}
	Let $\{x_{k}\}$ be the sequence generated by Algorithm \ref{Algorithmframework}. Then, we have 
	\begin{equation}\label{lemmaPhi}
		\begin{aligned}
			&\Phi_{k+1}(x_{k+1})-\Phi_{k}(x_{k})
			\leq-\min\left\{c_\lambda\|\rgrad\Phi_{k}(x_k)\|,{}c'\right\}\|\rgrad\Phi_{k}(x_k)\|+\nu_{k}+2\beta_kR^2,
		\end{aligned}
	\end{equation}
	{where the constant $c_\lambda$ is given after \eqref{sufficient decrease}.}
\end{lemma}
\begin{proof}{Proof}
	By the definition of $\Phi_k$ in \eqref{potential Phi_k} and the update formula for $y_{k+1}$ in \eqref{AlgorithmupdateY}, we have $\Phi_k(x_{k+1}) = F_k(x_{k+1}, y_{k+1})$. Similarly, by the definition of $y_{k + \frac{1}{2}}$ in \eqref{y_k+1/2}, we have $\Phi_{k+1}(x_{k+1}) = F_{k+1}(x_{k+1}, y_{k + \frac32})$, where we use $y_{k+\frac{3}{2}}$ as a shorthand for $y_{(k+1)+\frac{1}{2}}$. We now compute
	\begin{align}\label{suffcient decrease pf 1}
		{}\Phi_{k+1}(x_{k+1})-\Phi_{k}(x_{k+1})
		={}& F_{k+1}(x_{k+1},y_{k+\frac{3}{2}})-F_k(x_{k+1},y_{k+1})\nonumber\\
		\overset{(\text{a})}{\leq}{}& F_{k+1}(x_{k+1}, y_{k + \frac32}) - F_k(x_{k+1}, y_{k + \frac32})\nonumber \\
		\overset{(\text{b})}{=}{}&F(x_{k+1},y_{k+\frac{3}{2}})-\frac{\lambda}{2}\|y_{k+\frac{3}{2}}\|^2-\frac{\beta_{k+1}}{2}\|y_{k+\frac{3}{2}}-y_{k+1}\|^2\nonumber\\
		{}&-F(x_{k+1},y_{k+\frac{3}{2}})+\frac{\lambda}{2}\|y_{k+\frac{3}{2}}\|^2+\frac{\beta_{k}}{2}\|y_{k+\frac{3}{2}}-y_{k}\|^2\nonumber\\
		={}& \frac{\beta_{k}}{2}\|y_{k+\frac{3}{2}}-y_{k}\|^2-\frac{\beta_{k+1}}{2}\|y_{k+\frac{3}{2}}-y_{k+1}\|^2
		\overset{(\text{c})}{\leq}{} 2\beta_{k}R^2,
	\end{align}
	where (a) comes from the optimality of $y_{k+1}$ in \eqref{AlgorithmupdateY}, (b) holds by the definition of $F_{k+1}$ in \eqref{F_k}, and (c) is due to \eqref{R} and the fact that $y_k, y_{k+\frac32} \in \mathrm{dom}\,h$.
	By combining \eqref{suffcient decrease pf 1} and \eqref{sufficient decrease}, we obtain \eqref{lemmaPhi}.
\end{proof}

Based on Lemma \ref{lemma descent condition Phi}, we can establish the iteration complexity of Algorithm \ref{Algorithmframework} for returning an $\varepsilon$-RGS point of problem \eqref{prob:p1}.

\begin{theorem}\label{Theorem complexity minmax}
	Given a constant $\varepsilon > 0$, let  $\{(x_{k},y_{k})\} $ be the sequence generated by Algorithm \ref{Algorithmframework} with $ \lambda = {\varepsilon}/{(2R)}$. Suppose that Assumption \ref{assumption: level bound} holds. Then, there exists a $\hat{k}\leq K$ with 
	\begin{equation}\label{K Th1}
		K:=\left\lceil\left(\frac{4\beta_{1}R}{\varepsilon}\right)^{\frac{1}{\rho}}\right\rceil+ \left\lceil\frac{\Phi_1(x_1)-\Phi_*+{\Upsilon}}{\min\,\{c_\lambda\varepsilon,\,c'\}\,\varepsilon}\right\rceil
	\end{equation}
	such that $(x_{\hat{k}},y_{\hat{k}+\frac{1}{2}})$ is an $\varepsilon$-RGS point of problem \eqref{prob:p1}.
\end{theorem}
\begin{proof}{Proof}
	Given any $\gamma>0$, we first show that 
	$y_{k+\frac{1}{2}}$ satisfies 
	\begin{equation}\label{K_1}
		\frac{1}{\gamma}\left\| y_{k + \frac12} - \mathrm{prox}_{\gamma h}( y_{k + \frac12} + \gamma \mathcal{A}(x_k))\right\| \leq \varepsilon, \quad \forall\,k\geq K_1:=\left\lceil\left(\frac{4\beta_{1}R}{\varepsilon}\right)^{\frac{1}{\rho}}\right\rceil.
	\end{equation}
	For brevity, denote
	$
	y_{k+\frac{1}{2}}'=\prox_{\gamma h}(y_{k+\frac{1}{2}}+\gamma\mathcal{A}(x_{k}))
	$.
	It is easy to see that
	\begin{equation}\label{yk+1/2'}
		y_{k+\frac{1}{2}}' =\argmax_{y\in\mathcal{E}_2} \left\{F(x_{k},y)-\frac{1}{2\gamma}\|y-y_{k+\frac{1}{2}}\|^2\right\}.
	\end{equation}
	Due to the optimality of $y_{k+\frac12}$ in \eqref{y_k+1/2}, we have 
	\begin{equation}
		\begin{aligned}\label{lemma bound y-proxy 1}
			0
			\leq F_{k}(x_{k},y_{k+\frac{1}{2}})-F_{k}(x_{k},y_{k+\frac{1}{2}}' ).
		\end{aligned}
	\end{equation}
	Since $F(x_k,\cdot)-\frac{1}{2\gamma}\|\cdot-y_{k+\frac{1}{2}}\|^2$ is ${\gamma^{-1}}$-strongly concave, it follows from the optimality of $y_{k+\frac12}'$ in \eqref{yk+1/2'} that
	\begin{equation*}
		\begin{aligned}\label{lemma bound y-proxy 2}
			\frac{1}{2\gamma}\|y_{k+\frac{1}{2}}-y_{k+\frac{1}{2}}' \|^2
			\leq F(x_{k},y_{k+\frac{1}{2}}' )-\frac{1}{2\gamma}\|y_{k+\frac{1}{2}}' -y_{k+\frac{1}{2}}\|^2-F(x_{k},y_{k+\frac{1}{2}}). 
		\end{aligned}
	\end{equation*}
	This, together with \eqref{lemma bound y-proxy 1}, further implies that
	\begin{equation*}
		\begin{aligned}
			\frac{1}{\gamma}\|y_{k+\frac{1}{2}}-y_{k+\frac{1}{2}}' \|^2
			\leq&\, \left(F_k(x_k, y_{k + \frac12}) - F(x_k, y_{k + \frac12})\right) -  \left(F_k(x_k, y'_{k + \frac12}) - F(x_k, y'_{k + \frac12})\right)\\
			\overset{(\text{a})}{=}&\,\frac{\lambda}{2}\left(\|y_{k+\frac{1}{2}}'\|^2-\|y_{k+\frac{1}{2}}\|^2\right)+\frac{\beta_{k}}{2}\left(\|y_{k+\frac{1}{2}}'-y_{k}\|^2-\|y_{k+\frac{1}{2}}-y_{k}\|^2\right)\\
			=&\,\frac{\lambda}{2}\left\langle y_{k+\frac{1}{2}}'+y_{k+\frac{1}{2}},\,y_{k+\frac{1}{2}}'-y_{k+\frac{1}{2}}\right\rangle+\frac{\beta_{k}}{2}\left\langle y_{k+\frac{1}{2}}'+y_{k+\frac{1}{2}}-2y_{k},\,y_{k+\frac{1}{2}}'-y_{k+\frac{1}{2}}\right\rangle\\
			\leq&\,\left(\frac{\lambda}{2}\|y_{k+\frac{1}{2}}'+y_{k+\frac{1}{2}}\|+\frac{\beta_{k}}{2}\|y_{k+\frac{1}{2}}'+y_{k+\frac{1}{2}}-2y_{k}\|\right)\|y_{k+\frac{1}{2}}-y_{k+\frac{1}{2}}'\|,
		\end{aligned}
	\end{equation*}
	where (a) uses the definition of $F_k$ in \eqref{F_k}.
	Cancelling the term $\|y_{k + \frac12} - y'_{k + \frac12}\|$ on the two ends of the above inequalities yields
	\begin{equation}
			\frac{1}{\gamma}\|y_{k+\frac{1}{2}}-y_{k+\frac{1}{2}}' \|  \leq\frac{\lambda}{2}\|y_{k+\frac{1}{2}}'+y_{k+\frac{1}{2}}\|+\frac{\beta_{k}}{2}\|y_{k+\frac{1}{2}}'+y_{k+\frac{1}{2}}-2y_{k}\| 
 	\overset{(\text{a})}{\leq}\lambda R+2\beta_{k}R\overset{(\text{b})}{\leq}\frac{\varepsilon}{2}+\frac{\varepsilon}{2}=\varepsilon,
	\nonumber
	\end{equation}
	where (a) is due to \eqref{R} and the fact that $y_k,\,y_{k+\frac12},\,y_{k+\frac12}' \in \mathrm{dom}\,h$, and (b) uses $\lambda={\varepsilon}/{(2R)}$, $0\leq\beta_k \leq {\beta_1}/{k^\rho}$, and $k \geq K_1 = \big\lceil({4\beta_{1}R}/{\varepsilon})^{\frac{1}{\rho}}\big\rceil$.
	
	Next, we prove the existence of some $\hat{k} \in \{ K_1, K_1+1, \ldots, K \}$ such that $x_{\hat{k}}$ satisfies 
	\begin{equation}\label{Th1 proof1}
		\|\mathrm{grad}_x\, F(x_{\hat k}, y_{\hat k + \frac12})\| \leq \varepsilon. 
	\end{equation}
	Note that for any $x \in \mathcal{M}$, we deduce from \eqref{potential Phi_k} and the expression of $F(x,y)$ in \eqref{prob:p1} that 
	\begin{equation}\label{ineq: Phik geq Phi}
		\begin{aligned}
			\Phi_{k}(x)&=\max_{y\in\mathcal{E}_2}\left\{f(x)+\langle \mathcal{A}(x),y\rangle-h(y)-\frac{\lambda}{2}\|y\|^2-\frac{\beta_{k}}{2}\|y-y_k\|^2\right\}\\
			&\geq\max_{y\in\mathcal{E}_2}\left\{f(x)+\langle\mathcal{A}(x),y\rangle-h(y)-\frac{\lambda}{2}R^2-2\beta_{k}R^2\right\} \\
			&= \Phi(x)-\frac{\lambda}{2}R^2-2\beta_{k}R^2,
		\end{aligned}
	\end{equation}
	where the inequality is due to \eqref{R} and the last equality follows from the definition of $\Phi$ in \eqref{Phi}. 
	By \eqref{lemmaPhi}, we know that
	\begin{equation*}\label{lemma gradnormPhi}
		\begin{aligned}
			{}\min\left\{c_\lambda\|\rgrad\Phi_{k}(x_k)\|,\,c'\right\}\|\rgrad\Phi_{k}(x_k)\|
			\leq{}\Phi_{k}(x_{k})-\Phi_{k+1}(x_{k+1})+\nu_{k}+2\beta_{k}R^2.
		\end{aligned}
	\end{equation*}
	Summing the above inequality over $k= K_1, K_1 + 1, \ldots,K$ yields
	\begin{equation*}\label{Th1 proof sum}
		\begin{aligned}
			&\sum_{k=K_1}^{K}\min\left\{c_\lambda\|\rgrad\Phi_{k}(x_k)\|,\,c'\right\}\|\rgrad\Phi_{k}(x_k)\|\\
			\leq{}&\Phi_1(x_1)-\Phi_{K+1}(x_{K+1})+\sum_{k=1}^{K}(\nu_{k}+2\beta_{k}R^2)
			\leq\Phi_1(x_1)-\Phi^*+{\Upsilon},
		\end{aligned}
	\end{equation*}
	where the last inequality is due to Assumption \ref{assumption: level bound}, \eqref{Upsilon}, and \eqref{ineq: Phik geq Phi}.
	Hence, there exists some $\hat{k} \in \{K_1, K_1+1, \ldots, K \}$ such that
	\begin{equation}\label{Th1 proof bound gradPhi0}
		\min\left\{c_\lambda\|\rgrad\Phi_{\hat{k}}(x_{\hat{k}})\|,\,c'\right\}\|\rgrad\Phi_{\hat{k}}(x_{\hat{k}})\|\leq\frac{\Phi_1(x_1)-\Phi_*+{\Upsilon}}{K-K_1}.
	\end{equation}
	Substituting \eqref{K Th1} into \eqref{Th1 proof bound gradPhi0}, we have
	\begin{equation}\label{Th1 proof bound gradPhi}
		\|\rgrad\Phi_{\hat{k}}(x_{\hat{k}})\|\leq\varepsilon.
	\end{equation}
	The definition of $y_{k + \frac12}$ in \eqref{y_k+1/2}, together with \eqref{nabla Phik}, implies that $\nabla \Phi_{k}(x_{k}) = \nabla_x F(x_{k}, y_{k + \frac12})$ and thus 
	\begin{equation}\label{gradPhi=gradF}
		\rgrad\Phi_{k}(x_{k}) = \mathrm{grad}_x\,F(x_{k}, y_{k + \frac12}).
	\end{equation}
	Therefore, we see from \eqref{Th1 proof bound gradPhi} and \eqref{gradPhi=gradF} that \eqref{Th1 proof1} holds.
	Together with \eqref{K_1}, we conclude that $(x_{\hat k}, y_{\hat k + \frac12})$ is an $\varepsilon$-RGS point of problem \eqref{prob:p1}. 
\end{proof}

According to Theorem \ref{Theorem complexity minmax}, since $\lambda = {\varepsilon}/{(2R)}$ and $c_\lambda=\mathcal{O}(\lambda)$ by our choice, we see that the iteration complexity of Algorithm \ref{Algorithmframework} for obtaining an $\varepsilon$-RGS point of problem \eqref{prob:p1} is $\mathcal{O}(\varepsilon^{-3})$. 
It is worth mentioning that our iteration complexity result in Theorem \ref{Theorem complexity minmax} for the RADA algorithmic framework matches not only that for the less general ARPGDA in \cite{xu2023efficient2}, but also matches the best-known complexity results for algorithms that tackle NC-L minimax problems with $x$ being constrained on a nonempty compact convex set instead of a Riemannian manifold \cite{pan2021efficient,shen2023zeroth,he2024approximation}.
Furthermore, our iteration complexity analysis is based on the sufficient decrease condition \eqref{sufficient decrease} of the proposed value function $\Phi_{k}$. By contrast, the iteration complexity analyses in \cite{pan2021efficient,xu2023unified,xu2023efficient2} rely on the descent properties of the function $F(\cdot, y_k)$ and a different potential function, which do not extend to RADA. 

Next, we establish the iteration complexity of Algorithm \ref{Algorithmframework} for returning an $\varepsilon$-ROS point of problem \eqref{prob:p1}.
\begin{theorem}\label{Th2 complexity OS}
	Consider the setting of Theorem \ref{Theorem complexity minmax}. Then, $x_{\hat{k}}$ is an $\varepsilon$-ROS point of problem \eqref{prob:p1}, where $\hat{k}$ is the index given in Theorem \ref{Theorem complexity minmax}.
\end{theorem}
\begin{proof}{Proof}
	The proof is similar to that of Theorem \ref{Theorem complexity minmax}. 
	From \eqref{Th1 proof bound gradPhi}, we know there exists some $\hat{k} \in \{K_1, K_1+1, \ldots, K \}$ such that
	\begin{equation}\label{Th2 proof bound gradPhi}
		\|\rgrad\Phi_{\hat{k}}(x_{\hat{k}})\|\leq\varepsilon.
	\end{equation}
	Next, define $p_{k}=\prox_{(\lambda+\beta_{k}) h^*}(\mathcal{A}(x_{k})+\beta_{k}y_k)$. 
	From \eqref{gradient of Moreau envelope}, \eqref{nabla Phik},  and \eqref{y_k+1/2}, we have
	\begin{align}
		\rgrad\Phi_{k}(x_k)&=\rgrad f(x_k) + \proj_{\mathrm{T}_{x_k}\mathcal{M}}\left(\nabla \mathcal{A}(x_k)^\top y_{k+\frac12}\right)\quad\text{and}\quad
		y_{k+\frac12}\in \partial h^*(p_k).\label{equ: yk+12 in ph}
	\end{align}
	Then, it follows from \eqref{Th2 proof bound gradPhi} and \eqref{equ: yk+12 in ph} that
	\begin{equation}\label{Th2pf3}
		\begin{aligned}
			&\mathrm{dist}\left(0,\,\rgrad f(x_{\hat k}) + \proj_{\mathrm{T}_{x_{\hat{k}}}\mathcal{M}}\left(\nabla \mathcal{A}(x_{\hat k})^\top \partial h^*(p_{\hat{k}})\right)\right)\\
			\leq{}&\!\left\|\rgrad f(x_{\hat k}) + \mathrm{proj}_{\mathrm{T}_x \mathcal{M}}\left(\nabla\mathcal{A}(x_{\hat k})^\top y_{\hat{k}+\frac12}\right)\right\|
			=\left\|\rgrad \Phi_{\hat{k}}(x_{\hat{k}})\right\|
			\leq \varepsilon.
		\end{aligned}
	\end{equation}
	On the other hand, by the Moreau decomposition \eqref{Moreau Decomposition}, we have
	\begin{equation}\label{Th2pf4}
	\begin{aligned}
				\|p_{\hat{k}}-\mathcal{A}(x_{\hat{k}})\|=&\,\left\|\prox_{(\lambda+\beta_{\hat{k}})h^*}(\mathcal{A}(x_{\hat{k}})+\beta_{\hat{k}}y_{\hat{k}})-\mathcal{A}(x_{\hat{k}})\right\|\\
			\overset{}{=}&\,\left\|\beta_{\hat{k}}y_{\hat{k}}-(\lambda+\beta_{\hat{k}})\prox_{h/(\lambda+\beta_{\hat{k}})}\left(\frac{\mathcal{A}(x_{\hat{k}})+\beta_{\hat{k}}y_{\hat{k}}}{\lambda+\beta_{\hat{k}}}\right)\right\|\\
			\overset{(\text{a})}{\leq}&\,(\lambda+2\beta_{\hat{k}})R\overset{\text{(b)}}{\leq}\varepsilon,
	\end{aligned}
	\end{equation}
	where (a) uses $\mathrm{prox}_{h/(\lambda+\beta_{\hat{k}})} (\cdot) \in \mathrm{dom}\,h$ and \eqref{R}, and (b) uses $\lambda={\varepsilon}/{(2R)}$, $0\leq\beta_{\hat k} \leq {\beta_1}/{\hat k^\rho}$, and $\hat k \geq K_1 = \big\lceil({4\beta_{1}R}/{\varepsilon})^{\frac{1}{\rho}}\big\rceil$.
	Combining \eqref{Th2pf3} and \eqref{Th2pf4}, we conclude that $x_{\hat k}$ is an $\varepsilon$-ROS point of problem \eqref{prob:p1}. 
\end{proof}

Theorem \ref{Th2 complexity OS} shows that the iteration complexity of Algorithm \ref{Algorithmframework} for obtaining an $\varepsilon$-ROS point of problem \eqref{prob:p1} is $\mathcal{O}(\varepsilon^{-3})$ since $\lambda = {\varepsilon}/{(2R)}$ and $c_\lambda=\mathcal{O}(\lambda)$. 
This matches the best-known complexity results for the DSGM in \cite{beck2023dynamic} and the ManIAL in \cite{deng2024oracle}. However, the results in \cite{beck2023dynamic,deng2024oracle} apply only to the case where $\mathcal{A}$ is a linear mapping, whereas our results work for a general nonlinear mapping $\mathcal{A}$.

Before we leave this section, let us point out that by Theorems \ref{Theorem complexity minmax} and \ref{Th2 complexity OS}, Algorithm \ref{Algorithmframework} can find both an $\varepsilon$-RGS point and an $\varepsilon$-ROS point of problem \eqref{prob:p1} within $\mathcal{O}(\varepsilon^{-3})$ iterations. This is noteworthy, as the results in \cite{lin2020gradient,yang2022faster} for the Euclidean setting suggest that in general extra computation is needed to convert between an $\varepsilon$-RGS point and an $\varepsilon$-ROS point. Moreover, existing works on NC-L minimax problems \cite{pan2021efficient,he2024approximation,xu2023efficient2} only establish the iteration complexity for obtaining an $\varepsilon$-game-stationary point. The iteration complexity for obtaining an $\varepsilon$-optimization-stationary point is not addressed in these works.

\section{Two Efficient Algorithms for Problem \eqref{prob:p1}}\label{sec: two speicific}
In this section, we propose two efficient algorithms within the proposed RADA algorithmic framework, called RADA-PGD and RADA-RGD in Algorithms \ref{Algorithm ProjGrad} and \ref{Algorithm multiBB}, respectively, for solving the Riemannian NC-L minimax problem \eqref{prob:p1}. These two algorithms differ in the first-order method used for approximately solving the subproblem \eqref{subproblemphik} to find a point $x_{k+1} \in \mathcal{M}$ that satisfies the sufficient decrease condition \eqref{sufficient decrease}.

\begin{algorithm}[t]\label{Algorithm ProjGrad}
	\caption{RADA-PGD for problem \eqref{prob:p1}}
	
	Input $x_1\in\mathcal{M}, y_1\in\mathrm{dom}\,h$, $\lambda>0$, $\beta_{1}\geq0$, $\rho>1$, $\tau_1, \tau_2 \in (0,1)$, and a sequence $\{T_k\} \subset \mathbb{Z}_{+}$ with $1 \leq T_k \leq \overline{T}$, where $\overline{T}$ is a preset positive integer.
	
	\For{$k=1,\,2,\,\ldots$}
	{Set $x_{k,1}=x_{k}$.
		
		\For{$t=1,\,2,\,\ldots,\,T_k$}
		{Compute $\zeta_{k,t}=\ell_k^{-1}$ and update $x_{k,t+1}\in\proj_{\mathcal{M}}(x_{k,t}-\zeta_{k,t}\,\nabla\Phi_{k}(x_{k,t})).$
		}
		Update $x_{k+1}=x_{k,T_k+1}$, $y_{k+1}$ via \eqref{AlgorithmupdateY}, and $\beta_{k+1}$ via \eqref{RADA up beta}. }
\end{algorithm}
When the projection operator onto a closed manifold \(\mathcal{M} \) is tractable, which is the case in the motivating applications introduced in Section \ref{subsec: application}, we can employ the nonconvex projected gradient method to compute $x_{k+1}$ that satisfies the sufficient decrease condition \eqref{sufficient decrease}.
Noting that $\Phi_{k}$ is differentiable as shown in \eqref{nabla Phik}, it is straightforward to perform an ordinary projected gradient step with a suitable constant stepsize depending on the Lipschitz constant of $\nabla \Phi_k$. It is also possible to perform multiple projected gradient steps to achieve a more accurate solution to \eqref{subproblemphik} while satisfying \eqref{sufficient decrease}.
Within the RADA algorithmic framework, the first customized algorithm is RADA with projected gradient descent (RADA-PGD), which is formally presented in Algorithm \ref{Algorithm ProjGrad}. 
Here, $\mathbb{Z}_+$ is the set of all positive integers, and 
the integer $T_k \in [1, \overline{T}]$ denotes the number of the projected gradient descent steps to be performed at the $k$-th iteration with $\overline{T}$ being a predetermined positive integer.
The stepsize $\zeta_{k,t}$ is set as $\ell_k^{-1}$, where $\ell_k$ is to  be specified later in \eqref{ell_k}.

An alternative approach is to perform a Riemannian gradient descent step to obtain $x_{k+1}$ satisfying the sufficient decrease condition \eqref{sufficient decrease}.  For this method, determining a suitable constant stepsize is not as easy as in RADA-PGD.  Hence, for practical efficiency, we utilize the Barzilai-Borwein (BB) stepsize \cite{barzilai1988two,wen2013feasible,jiang2015framework,hu2018adaptive,iannazzo2018riemannian,jiang2022riemannian,gao2018new} as the initial stepsize and perform a backtracking line search to find a stepsize satisfying \eqref{sufficient decrease}.  
Similar to Algorithm \ref{Algorithm ProjGrad}, we can also perform multiple Riemannian gradient steps to achieve a more accurate solution to \eqref{subproblemphik}. 
The corresponding algorithm is named RADA with Riemannian gradient descent (RADA-RGD) and is presented in Algorithm \ref{Algorithm multiBB}, where the integer $T_k \in [1, \overline{T}]$ denotes the number of the Riemannian gradient descent steps to be performed at the $k$-th iteration and $\overline{T}$ is a predetermined positive integer. 
The Riemannian BB stepsize $\zeta_{k,t+1}^\mathrm{BB}$ in Algorithm \ref{Algorithm multiBB} is set as
\begin{equation}\label{BB stepsize}
	\zeta_{k,t+1}^\mathrm{BB}=
	\left\{
	\begin{aligned}
		\frac{\|x_{k,t+1}-x_{k,t}\|^2}{|\langle x_{k,t+1}-x_{k,t},\,v_{k,t}\rangle|}, &\quad \text { for odd } t; \\ 
		\frac{|\langle x_{k,t+1}-x_{k,t},\,v_{k,t}\rangle|}{\|v_{k,t}\|^2},&\quad \text { for even } t,
	\end{aligned}
	\right.
\end{equation}
where $v_{k,t}=\rgrad\Phi_{k}(x_{k,t+1})-\rgrad\Phi_{k}(x_{k,t})$.

Note that both RADA-PGD and RADA-RGD are single-loop algorithms, as the inner iteration number $T_k$ is bounded by a predetermined positive integer $\overline{T}$. The choice of $T_k$ is also flexible and can be adapted to the specific problem structure. Empirically, if the value function $\Phi_{k}$ is benign, for example, if the Lipschitz constant of its gradient can be well estimated, then setting $T_k\equiv1$ is often sufficient. Otherwise, we suggest setting a slightly larger $T_k$ to ensure a more substantial decrease in the value of $\Phi_{k}$, thereby improving the practical performance of RADA-PGD and RADA-RGD. This is particularly important because the subproblem \eqref{subproblemphik} is itself a minimax problem and can be viewed as a surrogate of problem \eqref{prob:p1}. Thanks to the proposed value function $\Phi_{k}$ and the flexible setting of $T_k$, our approaches to updating $x_{k+1}$ are expected to be more efficient than simply performing a gradient descent step on certain surrogate of $F(\cdot,y_k)$, as is done in \cite{xu2023unified,pan2021efficient,lu2020hybrid,he2024approximation,xu2023efficient2}.
This will be demonstrated in Section \ref{sec: numerical results}. Moreover, even when adapted to the Euclidean setting, our proposed RADA algorithmic framework and the associated RADA-PGD and RADA-RGD are new. 

\vspace{0.5cm}
\begin{algorithm}[t]\label{Algorithm multiBB}
	\caption{RADA-RGD for problem \eqref{prob:p1}}
	
	Input $x_1\in\mathcal{M}$, $y_1\in\mathrm{dom}\,h$, $\lambda > 0$, $\beta_{1}\geq0$, $0 < \zeta_{\min} < \zeta_{1,1} < \zeta_{\max}$, $\rho>1$, $\eta\in (0,1)$, $c_1\in(0,1)$, $\{\nu_{k}\}\in\mathbb{S}$, $\tau_1, \tau_2 \in (0,1)$, and a sequence $\{T_k\} \subset \mathbb{Z}_{+}$ with $1 \leq T_k \leq \overline{T}$, where $\overline{T}$ is a preset positive integer.
	
	\For{$k=1,\,2,\,\ldots$}
	{Set $x_{k,1}=x_{k}$. 
		
		\For{$t=1,\,2,\,\ldots,\,T_k$}
		{	Find the smallest nonnegative integer $j$ such that 
			\begin{equation}
				\begin{aligned}\label{linesearch condotion} 
					&\,\Phi_{k}\left(\mathrm{R}_{x_{k,t}}(-\zeta_{k,t}\eta^{j}\rgrad\,\Phi_{k}(x_{k,t}))\right)-\Phi_{k}(x_{k,t})
					\leq\, -c_1\zeta_{k,t}\eta^j\|\rgrad\,\Phi_{k}(x_{k,t})\|^2+\frac{\nu_k}{T_k}, 
				\end{aligned}
			\end{equation}
			and update $x_{k,t+1}=\mathrm{R}_{x_{k,t}}(-\zeta_{k,t}\eta^j\rgrad\,\Phi_{k}(x_{k,t})).$
			
			Compute  $\zeta_{k,t+1}^\mathrm{BB}$ according to \eqref{BB stepsize} and set
			\begin{equation}\label{zeta_kt}
			\zeta_{k,t+1}=\min\left\{\max\left\{\zeta_{k,t+1}^\mathrm{BB}, \zeta_{\min}\right\},  {\zeta_{\max}}/{\|\rgrad\Phi_{k}(x_{k,t+1})\|}\right\}.
		\end{equation}
		}
		{Update $x_{k+1}=x_{k,T_k+1}$, $y_{k+1}$ via \eqref{AlgorithmupdateY}, and $\beta_{k+1}$ via \eqref{RADA up beta}. } Set $\zeta_{k+1,1}= \zeta_{k,T_k+1}$.}
\end{algorithm}
\vspace{0.25cm}
\subsection{Verification of the Sufficient Decrease Condition \eqref{sufficient decrease}}
In this subsection, we show that the update $x_{k+1}$ in RADA-PGD and RADA-RGD satisfies the sufficient decrease condition \eqref{sufficient decrease}. Therefore, the iteration complexity of these two algorithms can be directly obtained from the results in Section \ref{sec: iteration complexity}. 

First, we make the following (somewhat standard) assumptions for our analysis; see, e.g., \cite{chen2024nonsmooth,chen2020proximal,boumal2019global,huang2022riemannian,huang2023inexact,liu2019quadratic,drusvyatskiy2019efficiency}.
\begin{assumption}\label{assumption: compact level set}
	The level set 
	\begin{equation*}\label{Omega}
		\Omega_{x_1} := \{x \in \mathcal{M} \mid \Phi(x) \leq \Phi(x_1)+\Upsilon\}
	\end{equation*}
	is compact, where $\Phi$ and $\Upsilon$ are defined in \eqref{Phi} and \eqref{Upsilon}, respectively, and $x_1$ is the initial point of Algorithm \ref{Algorithmframework}.
\end{assumption}
\begin{assumption}\label{assumption1}
	The function $f$ and the mapping $\mathcal{A}$ in \eqref{prob:p1} satisfy the following conditions:
	\begin{itemize}
		\item[(i)] The function $f: \mathcal{E}_1 \to \mathbb{R}$ is continuously differentiable and satisfies the descent property over $\mathcal{M}$, i.e., 
		\begin{equation}\label{descent Euclidean}
			f(x') \leq f(x) + \langle \nabla f(x),\,x' - x\rangle + \frac{L_f}{2} \|x' - x\|^2, \quad \forall\,x,x' \in \mathcal{M}. 
		\end{equation}
		
		\item[(ii)] The mapping $\mathcal{A}: \mathcal{E}_1 \to \mathcal{E}_2$ and its Jacobian mapping $\nabla \mathcal{A}$ are $L_{\mathcal{A}}^0$-Lipschitz and $L_{\mathcal{A}}^1$-Lipschitz over $\mathrm{conv}\,\mathcal{M}$, respectively. In other words, for any $x,x' \in \mathrm{conv}\,\mathcal{M}$, we have
		\begin{subequations}
			\begin{align}
				\| \mathcal{A}(x) - \mathcal{A}(x') \| \leq L_{\mathcal{A}}^0 \|x - x'\|,\label{L_A^0}\\
				\| \nabla \mathcal{A}(x) - \nabla \mathcal{A}(x') \| \leq L_{\mathcal{A}}^1 \|x - x'\|\label{L_A^1}. 
			\end{align}
		\end{subequations}
		Moreover, the mapping $\nabla \mathcal{A}$ is bounded over $\mathrm{conv}\,\mathcal{M}$, i.e., 
		\begin{equation}\label{rhoA}
			\rho_{\mathcal{A}}: = \sup_{x\,\in\,\mathrm{conv}\,\mathcal{M}}\, \|\nabla \mathcal{A} (x)\| < + \infty.
		\end{equation} 
	\end{itemize}
\end{assumption}
In our subsequent development, we allow for the possibility that the retraction at some $x \in \mathcal{M}$ is defined only locally; see \cite[Definition 2.1]{boumal2019global}. As such, we make the following assumption, which stipulates that the retraction is well defined on a compact subset of the tangent bundle; cf. \cite[Remark 2.2]{boumal2019global}.
\begin{assumption}\label{assmption: retraction well defined}
	There exists a constant $\varrho>0$ such that for any $(x,v)\in\mathcal{V}:=\{(x,v)\mid x\in\Omega_{x_1},\,v\in\mathrm{T}_x\mathcal{M},\,\|v\|\leq\varrho\}$, the point $\mathrm{R}_x(v)\in\mathcal{M}$ is well defined. 
\end{assumption}
The following lemma, which is extracted from \cite[Appendix B]{boumal2019global}, shows that the retraction satisfies the first- and second-order boundedness conditions on $\mathcal{V}$.
\begin{lemma}\label{lemma Bound retraction}
	Suppose that Assumptions \ref{assumption: compact level set} and \ref{assmption: retraction well defined} hold. Then, there exist constants  $\alpha_1, \alpha_2>0$ such that 
	\begin{equation*}\label{bound retraction}
		\|\mathrm{R}_x(v)-x\|\leq\alpha_1\|v\| \quad \text{{and}} \quad \|\mathrm{R}_x(v)-x-v\|\leq\alpha_2\|v\|^2
	\end{equation*} 
	for any $(x,v)\in\mathcal{V}$.
\end{lemma}

Based on Lemma \ref{lemma Bound retraction}, we can establish the descent property of $\Phi_k$ as follows. 
\begin{lemma}\label{lemma l-smooth E}
	Suppose that Assumptions \ref{assumption: compact level set}, \ref{assumption1}, and \ref{assmption: retraction well defined} hold. Then,
	the value function $\Phi_k$ in \eqref{potential Phi_k} satisfies the following properties: 
	\begin{itemize}
		\item[(i)] {\bf Euclidean Descent.} For any $x,x' \in \mathcal{M}$, we have
		\begin{equation}\label{lemma l-smooth E ineq E}
			{\Phi_k}(x') \leq {\Phi_k}(x) + \langle \nabla {\Phi_k}(x),\,x' - x \rangle + \frac{\ell_k}{2} \|x' - x\|^2, 
		\end{equation}
		where 
		\begin{equation}\label{ell_k}
			\ell_k= L_f+RL_{\mathcal{A}}^1+\frac{{\rho_\mathcal{A}}L_\mathcal{A}^0}{\lambda+\beta_{k}}.
		\end{equation}
		\item[(ii)] 
		{\bf Riemannian Descent.} For any $(x,v)\in\mathcal{V}$, we have
		\begin{equation}\label{lemma l-smooth E ineq R}
			{\Phi_k}(\mathrm{R}_x(v))\leq {\Phi_k}(x)+\langle \rgrad\, {\Phi_k}(x),\,v\rangle + \frac{L_k(x)}{2}\|v\|^2,
		\end{equation} 
		where 	
		\begin{equation} \label{equ:Lk} 
			L_k(x) = \ell_k \alpha_1^2 +2 \left(\|\nabla f(x)\|+\rho_\mathcal{A} R\right)\alpha_2.
		\end{equation}
	\end{itemize}
\end{lemma}
\begin{proof}{Proof}
	We first prove \eqref{lemma l-smooth E ineq E}. Let $P_k(x)=\Phi_{k}(x)-f(x)$.
	By \eqref{nabla Phik}, we have
	\begin{equation*}
		\nabla P_{k}(x)=\nabla\mathcal{A}(x)^\top\mathcal{B}(x),
	\end{equation*}
	where $$\mathcal{B}(x)=\prox_{{h}/(\lambda+\beta_{k})}\left(\frac{\mathcal{A}(x) + \beta_k y_k}{\lambda + \beta_k}\right).$$
	We claim that $\nabla P_{k}$ is Lipschitz continuous. Indeed, for any $x,x'\in\mathrm{conv}\,\mathcal{M}$, we have
	\begin{align*}
		\|\nabla P_{k}(x)-\nabla P_{k}(x')\|
		=\,&\| \nabla\mathcal{A}(x)^\top \mathcal{B}(x)- \nabla\mathcal{A}(x')^\top \mathcal{B}(x')\|\\
		\leq\,&\| (\nabla\mathcal{A}(x)-\nabla\mathcal{A}(x'))^\top \mathcal{B}(x)\|+\| \nabla\mathcal{A}(x')^\top (\mathcal{B}(x)-\mathcal{B}(x'))\|\\
		\overset{(\text{a})}{\leq}&\, L_\mathcal{A}^1\|\mathcal{B}(x)\|\cdot\|x-x'\|+\rho_\mathcal{A}\|\mathcal{B}(x)-\mathcal{B}(x')\|\\
		\overset{(\text{b})}{\leq}&\, L_\mathcal{A}^1\|\mathcal{B}(x)\|\cdot\|x-x'\|+\frac{\rho_\mathcal{A}}{\lambda+\beta_{k}}\|\mathcal{A}(x)-\mathcal{A}(x')\|\\
		\overset{(\text{c})}{\leq}&\, \left(RL_\mathcal{A}^1+\frac{\rho_\mathcal{A}L_\mathcal{A}^0}{\lambda+\beta_{k}}\right)\|x-x'\|,
	\end{align*}
	where (a) is due to \eqref{L_A^1} and \eqref{rhoA}, (b) is due to the nonexpansiveness of the proximal operator, and (c) is due to $\mathcal{B}(x) \in \mathrm{dom}\,h$, \eqref{R}, and \eqref{L_A^0}.
	Therefore, by \cite[Lemma 1.2.3]{nesterov2018lectures}, we know that $P_k$ satisfies the Euclidean descent property, namely, for any $x,x' \in \mathcal{M}$,
	\begin{equation}\label{descentpropert Pk}
		{P_k}(x') \leq {P_k}(x) + \langle \nabla {P_k}(x),\,x' - x \rangle + \frac{1}{2}\left(RL_{\mathcal{A}}^1+\frac{\rho_\mathcal{A}L_\mathcal{A}^0}{\lambda+\beta_{k}}\right) \|x' - x\|^2.
	\end{equation}
	Combining \eqref{descentpropert Pk} and \eqref{descent Euclidean} and noting that $\Phi_{k}=f+P_k$, we conclude that $\Phi_{k}$ satisfies the Euclidean descent property \eqref{lemma l-smooth E ineq E}.
	
	Moreover, using Lemma \ref{lemma Bound retraction} and following a similar analysis as that in \cite[Appendix B]{boumal2019global}, we know that $\Phi_k$ satisfies the Riemannian descent property \eqref{lemma l-smooth E ineq R}. This completes the proof.
\end{proof}

To prove that $x_{k+1}$ in Algorithm \ref{Algorithm ProjGrad} satisfies the sufficient decrease condition \eqref{sufficient decrease}, we need the following lemma.
\begin{lemma}\label{lemma PGD 1}
	Let $\left\{x_{k,t}\right\}_{t=1}^{T_k+1}$ be the sequence generated by Algorithm \ref{Algorithm ProjGrad} with $\zeta_{k,t}=\ell_k^{-1}$ at the $k$-th iteration. Suppose that 
	Assumptions \ref{assumption: compact level set}, \ref{assumption1}, and \ref{assmption: retraction well defined} hold. If $x_{k,t}\in\Omega_{x_1}$, then
	\begin{align*}
		\Phi_{k}(x_{k,t+1})-\Phi_{k}(x_{k,t})
		\leq-\min\left\{\frac{\|\rgrad\Phi_k(x_{k,t})\|}{2L_k(x_{k,t})},\,\frac{\varrho}{2}\right\}\|\rgrad\Phi_{k}(x_{k,t})\|.
	\end{align*}
\end{lemma}
\begin{proof}{Proof}
	By the definition of $x_{k,t+1}$ in Algorithm \ref{Algorithm ProjGrad}, we have 
	\begin{equation}\label{eq:xkt+1 is optimal}
		\begin{aligned}
			x_{k,t+1}\in{}&\argmin_{x\in\mathcal{M}}\left\|x-x_{k,t}+\frac{1}{\ell_k}\nabla\Phi_{k}(x_{k,t})\right\|^2 \\ 
			={}& \argmin_{x\in\mathcal{M}}\left\{\langle \nabla\Phi_{k}(x_{k,t}),\,x-x_{k,t}\rangle+\frac{\ell_k}{2}\|x-x_{k,t}\|^2\right\}.
		\end{aligned}
	\end{equation}
	Consider $x_{k,t+1}'=\mathrm{R}_{x_{k,t}}(d_{k,t})$ with $d_{k,t}=-c_{k,t}\,\rgrad\Phi_{k}(x_{k,t})$, where
	\begin{equation}\label{zetak}
		c_{k,t}:=\min\left\{\frac{1}{L_k(x_{k,t})},\,\frac{\varrho}{\|\rgrad\Phi_k(x_{k,t})\|}\right\}.
	\end{equation}
	We know that $(x_{k,t}, d_{k,t})\in\mathcal{V}$ from \eqref{zetak} where $\mathcal{V}$ is defined in Assumption \ref{assmption: retraction well defined}. By the optimality of $x_{k,t+1}$ in \eqref{eq:xkt+1 is optimal}, we get
	\begin{equation*} 
		\begin{aligned}
			&\langle \nabla\Phi_{k}(x_{k,t}),\,x_{k,t+1}-x_{k,t}\rangle+\frac{\ell_k}{2}\|x_{k,t+1}-x_{k,t}\|^2\\
			\leq{}&\langle \nabla\Phi_{k}(x_{k,t}),\,x_{k,t+1}'-x_{k,t}\rangle+\frac{\ell_k}{2}\|x_{k,t+1}'-x_{k,t} \|^2\\
			={}&\langle \nabla\Phi_{k}(x_{k,t}),\,d_{k,t}\rangle+\langle \nabla\Phi_{k}(x_{k,t}), x_{k,t+1}'-x_{k,t}-d_{k,t}\rangle+\frac{\ell_k}{2}\|x_{k,t+1}'-x_{k,t}\|^2\\
			\leq{}&\langle \nabla\Phi_{k}(x_{k,t}),\,d_{k,t}\rangle+\|\nabla\Phi_{k}(x_{k,t})\|\cdot\|x_{k,t+1}'-x_{k,t}-d_{k,t}\|+\frac{\ell_k}{2}\|x_{k,t+1}'-x_{k,t}\|^2\\
			\overset{(\text{a})}{\leq}& -c_{k,t}\|\rgrad\Phi_{k}(x_{k,t})\|^2 +\alpha_2 \left(\|\nabla f(x)\|+\rho_\mathcal{A} R\right)\|d_{k,t}\|^2+\frac{\alpha_1^2\ell_k}{2}\|d_{k,t}\|^2\\
			\overset{(\text{b})}{=}& 
			-c_{k,t} \left(1 - \frac{L_{k}(x_{k,t})}{2} c_{k,t}\right) \|\rgrad\Phi_{k}(x_{k,t})\|^2,
		\end{aligned}
	\end{equation*}
	where (a) follows from Lemma \ref{lemma Bound retraction}, $\mathrm{prox}_{h/(\lambda + \beta_k)}((\mathcal{A}(x) + \beta_k y_k)/(\lambda + \beta_k)) \in \mathrm{dom}\,h$, \eqref{nabla Phik}, \eqref{R}, and \eqref{rhoA}, and (b) is due to \eqref{equ:Lk}.
	In addition, applying \eqref{lemma l-smooth E ineq E} with $x' = x_{k,t+1}$ and $x =  x_{k,t}$ yields
	\begin{equation*}
		\begin{aligned}
			&\Phi_{k}(x_{k,t+1})-\Phi_{k}(x_{k,t})
			\leq\langle \nabla\Phi_{k}(x_{k,t}),\, x_{k,t+1}-x_{k,t}\rangle+\frac{\ell_k}{2}\|x_{k,t+1}-x_{k,t}\|^2\\
			\leq{}&	-c_{k,t}\left(1 - \frac{L_{k}(x_{k,t})}{2} c_{k,t}\right)\|\rgrad\Phi_{k}(x_{k,t})\|^2
			\overset{\text{(a)}}{\leq}-\frac{c_{k,t}}{2}\|\rgrad\Phi_{k}(x_{k,t})\|^2,
		\end{aligned}
	\end{equation*}
	where (a) is due to \eqref{zetak}. The desired result follows directly from substituting \eqref{zetak} into the above inequalities. 
\end{proof}

Based on Lemma \ref{lemma PGD 1}, we now present Proposition \ref{lemma RADA-PGD}, which shows that Algorithm \ref{Algorithm ProjGrad} is within the RADA algorithmic framework.
\begin{proposition}\label{lemma RADA-PGD}
	Let $\{x_{k}\}$ be the sequence generated by Algorithm \ref{Algorithm ProjGrad} with $\zeta_{k,t}=\ell_k^{-1}$. Suppose that Assumptions \ref{assumption: compact level set}, \ref{assumption1}, and \ref{assmption: retraction well defined} hold. Then, the sequence $\{x_k\}$ lies in the compact level set $\Omega_{x_1}$. Moreover, for $k=1,2,\ldots,$ we have
	\begin{equation}\label{suffcient decrease PGD}
		\Phi_{k}(x_{k+1})-\Phi_{k}(x_{k})\leq -\min\left\{\frac{\|\rgrad\Phi_k(x_{k})\|}{2L_k(x_k)},\,\frac{\varrho}{2}\right\}\|\rgrad\Phi_{k}(x_{k})\|.
	\end{equation} 
\end{proposition}
\begin{proof}{Proof}
	We apply mathematical induction. First, we show that \eqref{suffcient decrease PGD} holds when $k=1$. Since $x_{1,1}=x_1\in\Omega_{x_1}$, using Lemma \ref{lemma PGD 1} with $t=1$ and the definition of $c_{k,t}$ in \eqref{zetak}, we have
	\begin{align*}
		\Phi_{1}(x_{1,2})-\Phi_{1}(x_{1,1})
		\leq-\min\left\{\frac{\|\rgrad\Phi_1(x_{1,1})\|}{2L_1(x_{1,1})},\,\frac{\varrho}{2}\right\}\|\rgrad\Phi_{1}(x_{1,1})\|\leq0.
	\end{align*}
	This gives
	\begin{align}\label{x12 in level set}
		\Phi(x_{1,2})\overset{\text{(a)}}{\leq}\Phi_1(x_{1,2})+\left(\frac{\lambda}{2}+2\beta_{1}\right)R^2
		\overset{\text{(b)}}{\leq}\Phi_1(x_{1,1})+\Upsilon
		\overset{\text{(c)}}{\leq}\Phi(x_{1})+\Upsilon,
	\end{align}
	where (a) follows from \eqref{ineq: Phik geq Phi}, (b) uses the definition of $\Upsilon$ in \eqref{Upsilon}, and (c) holds by the definition of $\Phi$ in \eqref{Phi} and the fact that
	\begin{equation}\label{PhiklessthanPhi}
		\begin{aligned}
			\Phi_1(x)={}&\max_{y\in\mathcal{E}_2}\left\{f(x)+\langle \mathcal{A}(x),y\rangle-h(y)-\frac{\lambda}{2}\|y\|^2-\frac{\beta_{1}}{2}\|y-y_1\|^2\right\}\\
			\leq{}&\max_{y\in\mathcal{E}_2}\,\left\{f(x)+\langle \mathcal{A}(x),y\rangle-h(y)\right\}=\Phi(x).
		\end{aligned}
	\end{equation}
	Hence, we have $x_{1,2}\in\Omega_{x_1}$ from \eqref{x12 in level set}. Repeating the above procedure, we further know that
	\begin{align}\label{sufficient decrease Phi1 t}
		\Phi_{1}(x_{1,t+1})-\Phi_{1}(x_{1,t})
		\leq-\min\left\{\frac{\|\rgrad\Phi_1(x_{1,t})\|}{2L_1(x_{1,t})},\,\frac{\varrho}{2}\right\}\|\rgrad\Phi_{1}(x_{1,t})\|\leq0
	\end{align}
	and $x_{1,t +1} \in \Omega_{x_1}$ for $t=1,2,\ldots,T_1$.
	Summing the above inequality over $t=1,2,\ldots,T_1$ yields
	\begin{equation}\label{sufficient decrease Phi1}
		\begin{aligned}
			\Phi_{1}(x_{1,T_1+1})-\Phi_{1}(x_{1,1})={}&\sum_{t=1}^{T_1} \left(\Phi_1(x_{1,t+1})-\Phi_{1}(x_{1,t})\right) \\
			\leq{}&-\min\left\{\frac{\|\rgrad\Phi_1(x_{1})\|}{2L_1(x_1)},\,\frac{\varrho}{2}\right\}\|\rgrad\Phi_{1}(x_{1})\|.
		\end{aligned}
	\end{equation}
	Recalling that $x_{1,1} = x_1$ and $x_{1, T_1 + 1} = x_2$, it follows from \eqref{sufficient decrease Phi1} that \eqref{suffcient decrease PGD} holds for $k = 1$. 
	
	Next, we assume that $x_j \in \Omega_{x_1}$ for some $j \geq 2$ and \eqref{suffcient decrease PGD} holds for $k = 1,2,\ldots,j - 1$. It remains to prove that $x_{j+1} \in \Omega_{x_1}$ and \eqref{suffcient decrease PGD} holds for $k = j$. Combining \eqref{suffcient decrease pf 1} and \eqref{suffcient decrease PGD} with $k = 1, 2, \ldots, j -1$ yields 
	\begin{equation}\label{ineq: phik+1-phik}
		\Phi_{k+1}(x_{k+1})-\Phi_{k}(x_{k})\leq2\beta_{k}R^2, \quad \forall\,k=1,2,\ldots,j-1.
	\end{equation}
	Since $x_{j,1}=x_j\in\Omega_{x_1}$, Lemma \ref{lemma PGD 1} implies that
	\begin{equation}\label{inequ: Phij2}
		\Phi_j(x_{j,2})\leq\Phi_j(x_{j,1})=\Phi_j(x_j).
	\end{equation}
	Similar to the proof for $j=1$, we have
	\begin{equation}\label{inequ: Phij2 Upsilon}
		\begin{aligned}
			\Phi(x_{j,2})\leq{}&\Phi_{j}(x_{j,2})+\left(\frac{\lambda}{2}+2\beta_{j}\right)R^2\\
			\overset{\text{(a)}}{\leq}&\,\Phi_{j}(x_{j})+\left(\frac{\lambda}{2}+2\beta_{j}\right)R^2\\
			={}&\Phi_{1}(x_1)+\sum_{k=2}^{j}(\Phi_k(x_k)-\Phi_{k-1}(x_{k-1}))+\left(\frac{\lambda}{2}+2\beta_{j}\right)R^2\\
			\overset{(\text{b})}{\leq}&\, \Phi_1(x_1)+\frac{\lambda}{2}R^2+\sum_{k=1}^{+\infty}(\nu_{k}+2\beta_{k}R^2)\overset{(\text{c})}{\leq}\, \Phi(x_1)+\Upsilon,
		\end{aligned}
	\end{equation}
	where (a) is from \eqref{inequ: Phij2}, (b) is due to \eqref{ineq: phik+1-phik} and $\beta_j \leq \beta_1$, and (c) is due to \eqref{PhiklessthanPhi} and the definition of $\Upsilon$ in \eqref{Upsilon}. Hence, we get $x_{j,2}\in\Omega_{x_1}$. 
	Similar to \eqref{sufficient decrease Phi1 t} and \eqref{sufficient decrease Phi1}, we deduce from Lemma \ref{lemma PGD 1} that
	\begin{align}\label{suf de Phikt}
		\Phi_{j}(x_{j,t+1})-\Phi_{j}(x_{j,t})
		\leq-\min\left\{\frac{\|\rgrad\Phi_j(x_{j,t})\|}{2L_j(x_{j,t})},\,\frac{\varrho}{2}\right\}\|\rgrad\Phi_{j}(x_{j,t})\|\leq0
	\end{align}
	for $t=1,\,2,\,\ldots,\,T_j$. It follows from $x_{j,1} = x_j$ and $x_{j, T_j + 1} = x_{j+1}$  that
	\begin{equation*}
			\Phi_{j}(x_{j+1})-\Phi_{j}(x_j)=\sum_{t=1}^{T_j} \left(\Phi_j(x_{j,t+1})-\Phi_{j}(x_{j,t})\right) 
			\leq-\min\left\{\frac{\|\rgrad\Phi_j(x_{j})\|}{2L_j(x_j)},\,\frac{\varrho}{2}\right\}\|\rgrad\Phi_{j}(x_{j})\|.
	\end{equation*}
	Moreover, similar to \eqref{inequ: Phij2 Upsilon}, we have
	\begin{align*}
		\Phi(x_{j,T_j+1})\leq{}&\Phi_{j}(x_{j,T_j+1})+\left(\frac{\lambda}{2}+2\beta_{j}\right)R^2\\
		\overset{\text{(a)}}{\leq}{}&\Phi_{j}(x_{j})+\left(\frac{\lambda}{2}+2\beta_{j}\right)R^2\\
		={}&\Phi_{1}(x_1)+\sum_{k=2}^{j}(\Phi_k(x_k)-\Phi_{k-1}(x_{k-1}))+\left(\frac{\lambda}{2}+2\beta_{j}\right)R^2\\
		\overset{(\text{b})}{\leq}{}& \Phi_1(x_1)+\frac{\lambda}{2}R^2+\sum_{k=1}^{+\infty}(\nu_{k}+2\beta_{k}R^2)
		\overset{(\text{c})}{\leq}{} \Phi(x_1)+\Upsilon,
	\end{align*}
	where (a) is from \eqref{suf de Phikt} and $x_{j,1}=x_j$, (b) is due to \eqref{ineq: phik+1-phik} and $\beta_j \leq \beta_1$, and (c) is due to \eqref{PhiklessthanPhi} and the definition of $\Upsilon$ in \eqref{Upsilon}. Therefore, noting that $x_{j + 1} = x_{j, T_j + 1}$, we have $x_{j+1} \in \Omega_{x_1}$. Moreover,  similar to the derivation of \eqref{sufficient decrease Phi1}, we know from \eqref{suf de Phikt} that \eqref{suffcient decrease PGD} holds for $k = j$. This completes the inductive step and hence also the proof.
\end{proof}
\vspace{0.1cm}

Denote 
\begin{equation}
	\rho_\Phi:=\max_{x \in \Omega_{x_1}}\,\|\nabla f(x)\|+\rho_{\mathcal{A}}R\quad\text{and}\quad
	\underline{L}:=\ell_1\alpha_1^2+2\rho_{\mathcal{A}}R\alpha_2,\label{rhoPhi1}
\end{equation}
where $\rho_{\mathcal{A}}$ and $R$ are defined in \eqref{rhoA} and \eqref{R}, respectively. 
Then, for any $x\in\Omega_{x_1}$ and $k=1,\,2,\,\ldots,$ we have
\begin{equation}\label{rhoPhi ineq}
	\|\rgrad\Phi_{k}(x)\|\leq\|\nabla\Phi_{k}(x)\|\leq\rho_\Phi
\end{equation}
from \eqref{nabla Phik}, \eqref{R}, and \eqref{rhoA}. Moreover, we have
\begin{equation}\label{underbarLineq}
	0<\underline{L}\leq L_k(x)
\end{equation}
by substituting  $0\leq\beta_{k}\leq\beta_{1}$ into \eqref{equ:Lk}.
To prove that $x_{k+1}$ in Algorithm \ref{Algorithm multiBB} satisfies the sufficient decrease condition \eqref{sufficient decrease}, we need the following lemma.
\begin{lemma}\label{lemma RGD 1}
	Let $\left\{x_{k,t}\right\}_{t=1}^{T_k+1}$ be the sequence generated by Algorithm \ref{Algorithm multiBB} with $\zeta_{\max}\leq\varrho$ at the $k$-th iteration, where $\varrho$ is given in Assumption \ref{assmption: retraction well defined}. Suppose that 
	Assumptions \ref{assumption: compact level set}, \ref{assumption1}, and \ref{assmption: retraction well defined} hold. If $x_{k,t}\in\Omega_{x_1}$, then
	\begin{equation}\label{ineq: descent t RGD}
		\begin{aligned}
			\Phi_{k}(x_{k,t+1})-\Phi_{k}(x_{k,t})\leq - \frac{c_1 \min\left\{\rho_\Phi^{-1}  \zeta_{\max} \underline{L}, \zeta_{\min} \underline{L},  2\eta(1-c_1)\right\}}{L_k(x_{k,t})} \|\rgrad\Phi_{k}(x_{k,t})\|^2+\frac{\nu_k}{T_k}.
		\end{aligned}
	\end{equation}
\end{lemma}
\begin{proof}{Proof}
	Let $\hat{j}$ be the smallest nonnegative integer satisfying \eqref{linesearch condotion} and ${\hat\zeta_{k,t}}:=\zeta_{k,t}\eta^{\hat{j}}$. First, if $\hat j\geq1$, then by the backtracking line search procedure, we have
	\begin{equation}\label{ineq: bound hat zeta 1}
		\Phi_{k}(x_{k,t}^+)-\Phi_{k}(x_{k,t})>-c_1\frac{\hat\zeta_{k,t}}{\eta}\|\rgrad\Phi_{k}(x_{k,t})\|^2+\frac{\nu_{k}}{T_k}\geq -c_1\frac{\hat\zeta_{k,t}}{\eta}\|\rgrad\Phi_{k}(x_{k,t})\|^2,
	\end{equation}
	where $x_{k,t}^+:=\mathrm R_{x_{k,t}}\big(-{\hat\zeta_{k,t}}{\eta^{-1}}\rgrad\Phi_{k}(x_{k,t})\big)$.
	Since $\hat \zeta_{k,t} \eta^{-1}  = \zeta_{k,t} \eta^{\hat j - 1} \leq \zeta_{k,t}$ and \eqref{zeta_kt} implies that  $$\zeta_{k,t}\|\rgrad\Phi_{k}(x_{k,t})\|\leq\zeta_{\max}\leq\varrho,$$ we get $\big(x_{k,t},-{\hat\zeta_{k,t}}{\eta^{-1}}\rgrad\Phi_{k}(x_{k,t})\big) \in \mathcal{V}$, where $\mathcal{V}$ is defined in Assumption \ref{assmption: retraction well defined}.
	By Lemma \ref{lemma l-smooth E} (ii), we have 
	\begin{equation}\label{ineq: bound hat zeta 2}
		\Phi_{k}(x_{k,t}^+)-\Phi_{k}(x_{k,t})\leq-\left(1-\frac{\hat\zeta_{k,t}L_k(x_{k,t})}{2\eta}\right)\frac{\hat\zeta_{k,t}}{\eta}\|\rgrad\Phi_k(x_{k,t})\|^2.
	\end{equation}
	Combining \eqref{ineq: bound hat zeta 1} and \eqref{ineq: bound hat zeta 2} yields
	\begin{equation*}
		\hat\zeta_{k,t}\geq\frac{2\eta(1-c_1)}{L_k(x_{k,t})}.
	\end{equation*}
	Second, if $\hat j=0$, then $\hat\zeta_{k,t}=\zeta_{k,t}$. 
	We bound
	\begin{equation}\label{ineq bound lower zeta}
		\hat\zeta_{k,t}\geq\min\left\{\zeta_{k,t},\,\frac{2\eta(1-c_1)}{L_k(x_{k,t})}\right\}
		\geq \frac{\min\left\{\rho_\Phi^{-1} \zeta_{\max}  \underline{L}, \zeta_{\min} \underline{L},  2\eta(1-c_1)\right\}}{L_k(x_{k,t})},
	\end{equation}
	where the second inequality holds by \eqref{underbarLineq}, \eqref{rhoPhi ineq}, and $\zeta_{k,t} \geq \min\{\zeta_{\min}, \zeta_{\max}/\|\rgrad\Phi_{k}(x_{k,t+1})\|\}$ as defined in \eqref{zeta_kt}.
	By substituting \eqref{ineq bound lower zeta} into \eqref{linesearch condotion}, we get \eqref{ineq: descent t RGD}, as desired.
\end{proof}

Based on Lemma \ref{lemma RGD 1}, we have the following proposition, which shows that Algorithm \ref{Algorithm multiBB} is within the RADA algorithmic framework.
\begin{proposition}\label{lemma RADA-RGD}
	Let $\{x_{k}\}$ be the sequence generated by Algorithm \ref{Algorithm multiBB}. Suppose that Assumptions \ref{assumption: compact level set}, \ref{assumption1}, and \ref{assmption: retraction well defined}  hold. Then, the sequence $\{x_k\}$ lies in the compact level set $\Omega_{x_1}$. Moreover, for $k=1,2,\ldots,$ we have
	\begin{equation*}
		\Phi_{k}(x_{k+1})-\Phi_{k}(x_{k})\leq- \frac{c_1 \min\left\{\rho_\Phi^{-1}  \zeta_{\max} \underline{L}, \zeta_{\min} \underline{L},  2\eta(1-c_1)\right\}}{L_k(x_{k,t})}, \|\rgrad\Phi_{k}(x_{k})\|^2+\nu_k.
	\end{equation*} 
\end{proposition}
\begin{proof}{Proof}
	The result is obtained by using the same techniques as in the proof of Proposition \ref{lemma RADA-PGD}.
\end{proof}

Recalling that $\beta_k \geq 0$, we know from \eqref{ell_k}, \eqref{equ:Lk}, and \eqref{rhoPhi1} that 
\begin{equation*}\label{barL}
	L_k(x_k)\leq\bar{L}:=\alpha_1^2\left(L_f+RL_{\mathcal{A}}^1+\frac{{\rho_\mathcal{A}L_\mathcal{A}^0}}{\lambda}\right)+2\rho_\Phi\alpha_2=\mathcal{O}(\lambda^{-1}).
\end{equation*}
Therefore, from Propositions \ref{lemma RADA-PGD} and \ref{lemma RADA-RGD}, we see that  the sufficient decrease condition \eqref{sufficient decrease} holds with $$c_\lambda =\frac{1}{2\bar{L}}\quad\text{and}\quad c'=\frac{\varrho}{2}$$ for Algorithm \ref{Algorithm ProjGrad} and with 
\[
c_{\lambda}= \frac{\min\left\{\rho_\Phi^{-1} \zeta_{\max}  \underline{L}, \zeta_{\min} \underline{L},  2\eta(1-c_1)\right\}}{\bar L}, \quad\text{and}\quad c'= c_\lambda\rho_\Phi
\]  for Algorithm \ref{Algorithm multiBB}. Here, we use $\|\rgrad \Phi_k(x)\| = \min\,\{\|\rgrad \Phi_k(x)\|,\,\rho_{\Phi}\}$ as shown by \eqref{rhoPhi ineq}.
Consequently, by invoking Theorems \ref{Theorem complexity minmax} and \ref{Th2 complexity OS}, we obtain the following complexity result: 
\begin{theorem}
	Given a constant $\varepsilon > 0$, let  $\{(x_{k},y_{k})\} $ be the sequence generated by Algorithms \ref{Algorithm ProjGrad} or \ref{Algorithm multiBB} with $ \lambda = {\varepsilon}/{(2R)}$. Suppose that Assumptions \ref{assumption: level bound}, \ref{assumption: compact level set}, \ref{assumption1}, and \ref{assmption: retraction well defined} hold. Then, there exists an index $\hat{k}\leq K $ with $K = \mathcal{O}(\varepsilon^{-3})$ such that $(x_{\hat k}, y_{\hat {k}+\frac12})$ is an $\varepsilon$-RGS point and $x_{\hat{k}}$ is an $\varepsilon$-ROS point of problem \eqref{prob:p1}.
\end{theorem}

\section{Connections with Existing Algorithms}\label{sec: connections}
Recall that the minimax problem \eqref{prob:p1} is equivalent to the nonsmooth composite problem \eqref{prob:Riemannian nonsmooth composite}.
As it turns out, the proposed RADA algorithmic framework and the proposed RADA-RGD algorithm for solving problem \eqref{prob:p1} have interesting connections to the RALMs \cite{deng2023manifold,zhou2023semismooth,deng2024oracle} and RADMM \cite{li2023riemannian} for solving problem \eqref{prob:Riemannian nonsmooth composite}, respectively. In what follows, we elaborate on these connections and explain why RADA and RADA-RGD can be more advantageous than RALM and RADMM for tackling problem \eqref{prob:p1}.

\subsection{Connection between RADA and RALM}\label{Section relation RALM}

In this subsection, we show how the updates of the proposed RADA algorithmic framework for solving the minimax problem \eqref{prob:p1} relates to those of the RALMs \cite{deng2023manifold,zhou2023semismooth,deng2024oracle} for solving the equivalent nonsmooth composite problem \eqref{prob:Riemannian nonsmooth composite}. 

On one hand, using \eqref{Moreau Phi h}, the subproblem \eqref{subproblemphik} for obtaining the update $x_{k+1}$ in RADA can be written as 
\begin{equation}
	\min_{x \in \mathcal{M}} \left\{f(x) + M_{{(\lambda+\beta_{k})h^*}}(\mathcal{A}(x)+\beta_{k}y_{k})\right\}.\label{Moreau Phi hstar}
\end{equation}
Moreover, by the Moreau decomposition \eqref{Moreau Decomposition}, the update formula \eqref{AlgorithmupdateY} for $y_{k+1}$ in RADA can be written as  
\begin{equation}\label{Alg1upYprox}
	\begin{aligned}
		y_{k+1} 
		&=\frac{1}{\lambda+\beta_{k}}\left(\beta_{k}y_{k}+\mathcal{A}(x_{k+1})-\prox_{(\lambda+\beta_{k})h^*}(\mathcal{A}(x_{k+1})+\beta_{k}y_{k}) \right).
	\end{aligned}
\end{equation}
On the other hand, consider applying the RALM in \cite{deng2023manifold,zhou2023semismooth,deng2024oracle} to the following equivalent reformulation of \eqref{prob:Riemannian nonsmooth composite}:
\begin{equation}\label{prob: eq prob for ALM}
	\min_{x\in\mathcal{M},\,q\in\mathcal{E}_2}\left\{f(x) +h^*(q)\right\} \quad \mathrm{s.t.} \quad \mathcal{A}(x) - q = 0.
\end{equation}
Here, $q$ is an auxiliary variable.
The augmented Lagrangian function associated with problem \eqref{prob: eq prob for ALM} is given by
\begin{equation*}\label{augmented Lagrangian}
	\mathcal{L}_\sigma(x, q;y):= f(x)+h^*(q)+\langle y,\,\mathcal{A}(x)- q\rangle+\frac{\sigma}{2}\|\mathcal{A}(x)- q\|^2,
\end{equation*}
where $y$ is the Lagrange multiplier corresponding to the constraint $\mathcal{A}(x) -q = 0$ and $\sigma > 0$ is the penalty parameter.
At the $k$-th iteration, instead of utilizing the subproblem
\begin{equation}\label{prob: RALM subprob}
	{\min_{x \in \mathcal{M},\,q\in\mathcal{E}_2} \mathcal{L}_{\sigma_k}(x, q;y_k),}
\end{equation}
RALM  utilizes the equivalent subproblem
\begin{equation}\label{RALM subproblem}
	\min_{x\in\mathcal{M}}\left\{\mathcal{L}_k(x):=f(x)+M_{{h^*/\sigma_k}}\left(\mathcal{A}(x)+\frac{{y_k}}{\sigma_{k}}\right)\right\}
\end{equation}
to obtain the update $x_{k+1}$, where the equivalence follows by eliminating the variable $q$ in \eqref{prob: RALM subprob}. Specifically, RALM first approximately solves \eqref{RALM subproblem} to obtain a point $x_{k+1}\in\mathcal{M}$ that satisfies
\begin{equation}\label{equ:RALM:x}
	\|\rgrad\mathcal{L}_k(x_{k+1})\|\leq\varepsilon_k,
\end{equation}
where the sequence of tolerances $\left\{ \varepsilon_k \right\}$ converges to 0. 
Then, it updates
\begin{equation}\label{RALMupY}
	\begin{aligned}
		y_{k+1}= y_k+\sigma_k\left(\mathcal{A}(x_{k+1})-\prox_{{{h^*}/{\sigma_k}}}\left(\mathcal{A}(x_{k+1})+\frac{y_k}{\sigma_{k}}\right)\right).
	\end{aligned}
\end{equation}

Now, let us compare the subproblems \eqref{Moreau Phi hstar} and \eqref{RALM subproblem} for updating $x$ and the formulas \eqref{Alg1upYprox} and \eqref{RALMupY} for updating $y$.
When $\lambda = 0$ and $\beta_k = {\sigma_k^{-1}}$, 
the subproblem \eqref{Moreau Phi hstar} in RADA coincides with the subproblem \eqref{RALM subproblem} in RALM, and the update formula \eqref{Alg1upYprox} in RADA reduces to the update formula \eqref{RALMupY} in RALM. 
However, when $\lambda > 0$, both \eqref{Moreau Phi hstar} and \eqref{Alg1upYprox} in RADA differ from their respective counterparts \eqref{RALM subproblem} and \eqref{RALMupY} in RALM. More importantly, RADA and RALM differ in their stopping criteria for the subproblems. 
Indeed, RADA only needs to find a solution $x_{k+1}$ that satisfies the sufficient decrease condition \eqref{sufficient decrease}, whereas RALM needs to find a solution $x_{k+1}$ that satisfies the inexactness condition \eqref{equ:RALM:x}. The former offers significantly more flexibility when solving the subproblem.
It should be mentioned that \eqref{Moreau Phi hstar} is also related to the so-called dampened or perturbed augmented Lagrangian function (see, e.g., in \cite{hajinezhad2019perturbed,kong2024global}), which is used in the design of primal-dual algorithms for linearly constrained composite optimization problems in the Euclidean space.
Furthermore, when $\beta_{k}\equiv0$ and $\lambda>0$, the subproblem \eqref{Moreau Phi hstar} becomes
\begin{equation}\label{prob: nonsmooth comp 2}
	\min_{x \in \mathcal{M}}\left\{f(x)+M_{\lambda h^*}(\mathcal{A}(x))\right\},
\end{equation}
which corresponds to the subproblem in the Riemannian smoothing methods in \cite{peng2023riemannian,beck2023dynamic,zhang2023riemannian}. 

The above relationship between the proposed RADA algorithmic framework and RALM provides important insights into the update formula \eqref{RADA up beta} for $\beta_{k+1}$ based on the parameter $\delta_{k+1}$ defined in \eqref{delta}, as well as the numerical performance of the two approaches. Specifically, at the $k$-th iteration of RALM, we can express the update $q_{k+1}$ as
\begin{equation*}
	q_{k+1}=\prox_{{{h^*}/{\sigma_k}}}\left(\mathcal{A}(x_{k+1})+\frac{y_k}{\sigma_{k}}\right).
\end{equation*}
Using \eqref{Alg1upYprox}, we see that our Algorithm \ref{Algorithmframework} has an analogous update, namely,
\begin{equation*}
	p_{k+1}=\prox_{(\lambda+\beta_{k})h^*}\left(\mathcal{A}(x_{k+1})+\beta_{k}y_k\right).
\end{equation*}
Now, using \eqref{delta}, we compute
\begin{equation*}\label{primal residual}
	\begin{aligned}
		\delta_{k+1}=\,&\|(\lambda+\beta_{k})y_{k+1}-\beta_{k}y_k\|_{\infty}\\
		\overset{(\text{a})}{=}&\left\|-\beta_{k}y_k+(\lambda+\beta_{k})\prox_{{h}/(\lambda+\beta_{k})}\left(\frac{\mathcal{A}(x_{k+1})+\beta_{k}y_k}{\lambda+\beta_{k}}\right)\right\|_{\infty}\\
		\overset{(\text{b})}{=}&\|\mathcal{A}(x_{k+1})-\prox_{(\lambda+\beta_{k})h^*}(\mathcal{A}(x_{k+1})+\beta_{k}y_k)\|_{\infty}\\
		=\,&\|\mathcal{A}(x_{k+1})-p_{k+1}\|_{\infty},
	\end{aligned}
\end{equation*}
where (a) is due to \eqref{AlgorithmupdateY} and (b) follows from the Moreau decomposition \eqref{Moreau Decomposition}. This shows that the parameter $\delta_{k+1}$ in \eqref{delta} can be viewed as a primal residual in RALM.

\subsection{Connection between RADA-RGD and RADMM}\label{subsection: Interpreting RADA-RGD from the Perspective of ADMM}

In this subsection, 
we first elaborate on the relationship between RADA-RGD with $T_k \equiv 1$ and the RADMM in \cite{li2023riemannian}. We then show that RADA-RGD with $T_k \equiv 1$ is equivalent to a new Riemannian sGS-ADMM. 

Consider applying the RADMM in \cite{li2023riemannian} to the following equivalent formulation of the smoothed problem \eqref{prob: nonsmooth comp 2}:
\begin{equation}\label{prob: smoothed nonsmoooth composite}
	\min_{x\in\mathcal{M}}\left\{f(x) +M_{{\lambda h^*}}( p)\right\} \quad \mathrm{s.t.} \quad \mathcal{A}(x) - p = 0.
\end{equation}
The augmented Lagrangian function associated with problem \eqref{prob: smoothed nonsmoooth composite} is given by
\begin{equation*}
	\tilde{\mathcal{L}}_\sigma(x,p;y):= f(x)+M_{{\lambda h^*}}(p)+\langle y,\mathcal{A}(x)-p\rangle+\frac{\sigma}{2}\|\mathcal{A}(x)-p\|^2,
\end{equation*}
where $y$ is the Lagrange multiplier corresponding to the constraint $\mathcal{A}(x) - p = 0$ and $\sigma > 0$ is the penalty parameter.
At the $k$-th iteration, RADMM computes the updates
\begin{subequations}\label{RADMM}
	\begin{align}
		x_{k+1} = &\, \mathrm{R}_{x_{k}} (-\zeta_{k}\,\mathrm{grad}_x\,\tilde{\mathcal{L}}_{\sigma_{k}}(x_k,p_k;y_k)), \label{RADMMx} \\
		p_{k+1}= &\, \argmin_{p\in\mathcal{E}_2}\,  \tilde{\mathcal{L}}_{\sigma_{k}}(x_{k+1},p;y_k), \\
		y_{k+1} = &\, y_{k}+\sigma_{k}(\mathcal{A}(x_{k+1})-p_{k+1}),
	\end{align}
\end{subequations}
where the stepsize $\zeta_{k}>0$ in \eqref{RADMMx} and the penalty parameter $\sigma_{k}>0$ are both taken to be constants in \cite{li2023riemannian}. 

Next, we modify RADA-RGD in Algorithm \ref{Algorithm multiBB} by introducing an additional variable $p$, so that the resulting algorithm, which is shown in Algorithm \ref{Algorithm multiBB2}, can be used to solve problem \eqref{prob:Riemannian nonsmooth composite}. 
Here, for simplicity of presentation, we omit the line search procedure in RADA-RGD.

\begin{algorithm}[t]\label{Algorithm multiBB2}
	\caption{RADA-RGD reformulation for problem \eqref{prob:Riemannian nonsmooth composite}}
	
	Input $x_1\in\mathcal{M}$, $y_1\in\mathcal{E}_2$, $\lambda > 0$, $\tau_1, \tau_2 \in (0,1)$, and a sequence $\{T_k\} \subset \mathbb{Z}_{+}$ with $1 \leq T_k \leq \overline{T}$, where $\overline{T}$ is a preset positive integer.
	
	\For{$k=1,\,2,\,\ldots$}
	{Set $x_{k,1}=x_{k}$. 
		
		\For{$t=1,\,2,\,\ldots,\,T_k$}
		{
			Compute $\zeta_{k,t}$ as in Algorithm \ref{Algorithm multiBB} and update 
			\begin{align}
				p_{k,t}&=\argmin_{p\in\mathcal{E}_2}\,  \tilde{\mathcal{L}}_{1/\beta_{k}}(x_{k,t},p;y_k),\nonumber\\
				x_{k,t+1}&=\mathrm{R}_{x_{k,t}}(-\zeta_{k,t}\,\mathrm{grad}_x\,\tilde{\mathcal{L}}_{1/\beta_{k}}(x_{k,t},p_{k,t};y_k)).\label{RGDupx2}
			\end{align}
		}
		{Set $x_{k+1}=x_{k,T_k+1}$.
			
			Update $p_{k+1}=\argmin_{p\in\mathcal{E}_2}\,  \tilde{\mathcal{L}}_{1/\beta_{k}}(x_{k+1},p;y_k)$
			and	\begin{equation}\label{RGDupy2}
				y_{k+1}=y_k+\frac{1}{\beta_k}(\mathcal{A}(x_{k+1})-p_{k+1}).
			\end{equation}
			
			Update $\beta_{k+1}$ as \eqref{RADA up beta}.
		} 
	}
\end{algorithm}
The key difference between \eqref{RADMM} and Algorithm \ref{Algorithm multiBB2} with $T_k \equiv 1$ is that the latter updates $x$ and $p$ in a symmetric Gauss-Seidel fashion, i.e., $p \to x \to p$.
Therefore, Algorithm \ref{Algorithm multiBB2}
with $T_k \equiv 1$ can be viewed as a Riemannian sGS-ADMM (cf. \cite{chen2017efficient,li2016schur,li2019block} for the sGS-ADMM in Euclidean convex settings) for solving the Riemannian nonsmooth composite problem \eqref{prob:Riemannian nonsmooth composite}.  
To the best of our knowledge, this is the first time a Riemannian sGS-type ADMM with a convergence guarantee is proposed.

Now, let us establish the equivalence between the original RADA-RGD (i.e., Algorithm \ref{Algorithm multiBB}) and its modified version (i.e., Algorithm \ref{Algorithm multiBB2}).

\begin{proposition}\label{relation ADMM}
	The updates of $x$ and $y$ in Algorithm \ref{Algorithm multiBB} are equivalent to \eqref{RGDupx2} and \eqref{RGDupy2} in Algorithm \ref{Algorithm multiBB2}, respectively. 
\end{proposition}
\begin{proof}{Proof}
	Denote $p(x):=\argmin_{p\in\mathcal{E}_2}\,  \tilde{\mathcal{L}}_{1/\beta_{k}}(x,p;y_k)$.
	By \cite[Lemma 1]{li2023riemannian}, we have
	\begin{align*}
		p(x) & = \frac{\beta_{k}}{\lambda+\beta_{k}}\prox_{(\lambda+\beta_{k})h^*}\left(\mathcal{A}(x)+\beta_{k}y_{k}\right) + \frac{\lambda}{\lambda+\beta_{k}}(\mathcal{A}(x)+\beta_{k}y_{k}), \label{PinNADMM}
	\end{align*}
	which implies that
	\begin{equation}\label{A-p}
		\mathcal{A}(x)-p(x)=\frac{\beta_{k}}{\lambda+\beta_{k}}\left(\mathcal{A}(x)-\lambda y_k-\prox_{(\lambda+\beta_{k})h^*}\left(\mathcal{A}(x)+\beta_{k}y_{k}\right) \right).
	\end{equation}
	We first show the equivalence of the updates of $y$ in the two algorithms. 
	Note that $p(x_{k+1})=p_{k+1}$.
	Substituting \eqref{A-p} with $x = x_{k+1}$ into \eqref{RGDupy2} yields
	\begin{equation*}\label{NADMMupdateYsmooth}
		\begin{aligned}
			y_{k+1}  
			=&\, \frac{\beta_{k}}{\lambda+\beta_{k}}y_{k}+\frac{1}{\lambda+\beta_{k}}\left(\mathcal{A}(x_{k+1})-\prox_{(\lambda+\beta_{k})h^*}(\mathcal{A}(x_{k+1})+\beta_{k}y_{k})\right),
		\end{aligned}
	\end{equation*}
	which is exactly the update \eqref{Alg1upYprox}, or equivalently, the update \eqref{AlgorithmupdateY} in Algorithm \ref{Algorithm multiBB}. 
	
	We now show the equivalence of the updates of $x$ in the two algorithms. Substituting \eqref{A-p} into $\nabla_x \tilde{\mathcal{L}}_{1/\beta_{k}}(x,p(x);y_k)$ yields
	\begin{equation*}\label{nablaL}
		\begin{aligned}
			 \nabla_x \tilde{\mathcal{L}}_{1/\beta_{k}}(x,p(x);y_k) = &\, \nabla f(x)+\nabla \mathcal{A}(x)^\top\left(y_{k}+\frac{1}{\beta_{k}}(\mathcal{A}(x)-p(x))\right)\\ 
			=&\, \nabla f(x)+\frac{1}{\lambda+\beta_{k}}\nabla \mathcal{A}(x)^\top\left(\mathcal{A}(x)+\beta_{k}y_{k}-\prox_{(\lambda+\beta_{k})h^*}(\mathcal{A}(x)+\beta_{k}y_{k})\right),
		\end{aligned}
	\end{equation*}
	which, together with \eqref{nabla Phik} and the Moreau decomposition \eqref{Moreau Decomposition}, implies that 
	\[
	\nabla_x \tilde{\mathcal{L}}_{1/\beta_{k}}(x,p(x);y_k) = \nabla \Phi_k(x). 
	\]
	Therefore, noting that $p(x_{k,t}) = p_{k,t}$, we can rewrite \eqref{RGDupx2} as 
	\[
	x_{k,t+1}=\mathrm{R}_{x_{k,t}}(-\zeta_{k,t}\,\rgrad\Phi_{k}(x_{k,t})),
	\]
	which is exactly the update of $x$ in RADA-RGD. 
\end{proof}

\section{Numerical Results}\label{sec: numerical results} 
In this section, we report the numerical results of our proposed RADA algorithmic framework for solving SPCA, FPCA, and SSC problems. 
All the tests are implemented in MATLAB 2023b and evaluated on Apple M2 Pro CPU. Our code is available at \url{https://github.com/XuMeng00124/RADAopt}.

\begin{table}
	\centering
	\fontsize{10pt}{\baselineskip}\selectfont
	\caption{Average performance comparison on SPCA.}
	\tabcolsep=0.15cm
	\renewcommand\arraystretch{1}
	\begin{tabular}{cccccc|ccccc|ccccc}
		\hline
		&\multicolumn{5}{c}{RADA-RGD}  & \multicolumn{5}{c}{RALM-II} & \multicolumn{5}{c}{ManPG-Ada}  \\
		\cline{2-16}
		& $-\Phi$ & var & iter & cpu  & spar  & $-\Phi$ & var & iter & cpu & spar  & $-\Phi$ & var & iter & cpu & spar  \\
		\hline
		$\mu$ & \multicolumn{15}{c}{Synthetic, $d=1000, N=50, r=10$}\\
		\hline
		0.5 & $\mathbf{859.3}$ & $0.991$ & $636$ & $\mathbf{1.9}$ & $26.6$ &$859.0$ & $0.991$ & $26$ & $5.4$ & $26.7$ &$859.1$ & $0.991$ & $18366$ & $11.1$ & $26.5$ \\ 
		0.75 & $\mathbf{810.2}$ & $0.979$ & $2288$ & $\mathbf{1.3}$ & $35.5$ &$809.9$ & $0.978$ & $28$ & $5.3$ & $35.7$ &$810.1$ & $0.978$ & $19581$ & $12.0$ & $35.5$ \\ 
		1 & $\mathbf{765.8}$ & $0.965$ & $494$ & $\mathbf{1.5}$ & $42.1$ &$765.5$ & $0.964$ & $30$ & $5.2$ & $42.3$ &$765.4$ & $0.964$ & $24640$ & $14.5$ & $42.3$ \\ 
		1.25 & $\mathbf{724.2}$ & $0.952$ & $569$ & $\mathbf{1.7}$ & $47.4$ &$723.8$ & $0.951$ & $33$ & $5.3$ & $47.4$ &$724.0$ & $0.951$ & $18271$ & $12.0$ & $47.4$ \\ 
		1.5 & $\mathbf{685.8}$ & $0.936$ & $527$ & $\mathbf{1.6}$ & $51.8$ &$685.7$ & $0.936$ & $36$ & $5.4$ & $51.9$ &$685.3$ & $0.936$ & $13713$ & $9.1$ & $51.9$ \\ 
		\hline
		$r$ & \multicolumn{15}{c}{Synthetic, $d=2000, N=50, \mu=3$}\\
		\hline        
		8 & $\mathbf{1163.1}$ & $0.919$ & $532$ & $\mathbf{2.5}$ & $46.1$ &$1161.6$ & $0.923$ & $61$ & $16.2$ & $47.6$ &$1161.0$ & $0.915$ & $15004$ & $12.4$ & $49.2$ \\ 
		10 & $\mathbf{1195.2}$ & $0.904$ & $624$ & $\mathbf{3.7}$ & $59.3$ &$1193.6$ & $0.896$ & $62$ & $21.1$ & $59.8$ &$1190.0$ & $0.885$ & $16399$ & $19.2$ & $61.0$ \\ 
		12 & $\mathbf{1189.7}$ & $0.890$ & $816$ & $\mathbf{5.5}$ & $65.1$ &$1188.7$ & $0.891$ & $70$ & $29.6$ & $64.3$ &$1188.6$ & $0.890$ & $12614$ & $22.6$ & $65.2$ \\ 
		14 & $\mathbf{1145.9}$ & $0.867$ & $698$ & $\mathbf{5.7}$ & $69.8$ &$1144.5$ & $0.869$ & $64$ & $31.5$ & $69.0$ &$1144.6$ & $0.865$ & $15179$ & $46.8$ & $69.9$ \\ 
		16 & $\mathbf{1247.8}$ & $0.895$ & $739$ & $\mathbf{7.0}$ & $75.6$ &$1245.0$ & $0.892$ & $50$ & $28.5$ & $76.1$ &$1240.1$ & $0.894$ & $23799$ & $98.3$ & $76.0$ \\ 
		\hline
		$\mu$ & \multicolumn{15}{c}{coil-100, $d=1024, N=100, r=10$}\\
		\hline
		1 & $\mathbf{531.0}$ & $0.960$ & $268$ & $\mathbf{1.1}$ & $47.2$ &$530.8$ & $0.959$ & $34$ & $5.6$ & $47.4$ &$530.4$ & $0.959$ & $17205$ & $11.9$ & $47.3$ \\ 
		1.25 & $\mathbf{488.5}$ & $0.940$ & $234$ & $\mathbf{1.0}$ & $55.3$ &$487.8$ & $0.940$ & $35$ & $4.9$ & $55.2$ &$487.6$ & $0.940$ & $11865$ & $8.9$ & $55.2$ \\ 
		1.5 & $\mathbf{453.5}$ & $0.923$ & $266$ & $\mathbf{1.2}$ & $61.3$ &$451.1$ & $0.917$ & $39$ & $5.1$ & $61.7$ &$448.5$ & $0.911$ & $11828$ & $8.7$ & $61.7$ \\ 
		1.75 & $\mathbf{419.7}$ & $0.906$ & $319$ & $\mathbf{1.4}$ & $66.4$ &$417.2$ & $0.898$ & $41$ & $5.7$ & $66.5$ &$411.5$ & $0.885$ & $10475$ & $7.7$ & $66.5$ \\ 
		2 & $\mathbf{388.2}$ & $0.890$ & $344$ & $\mathbf{1.5}$ & $70.6$ &$380.5$ & $0.866$ & $43$ & $5.6$ & $70.2$ &$373.7$ & $0.848$ & $10950$ & $7.9$ & $70.4$ \\ 
		\hline
		$r$ & \multicolumn{15}{c}{coil-100, $d=1024, N=100, \mu=2$}\\
		\hline        
		8 & $\mathbf{384.8}$ & $0.913$ & $254$ & $\mathbf{0.9}$ & $63.2$ &$380.3$ & $0.898$ & $40$ & $3.9$ & $63.5$ &$375.2$ & $0.883$ & $16313$ & $7.1$ & $64.2$ \\ 
		10 & $\mathbf{388.2}$ & $0.890$ & $344$ & $\mathbf{1.5}$ & $70.6$ &$380.5$ & $0.866$ & $43$ & $5.7$ & $70.2$ &$373.7$ & $0.848$ & $10950$ & $7.9$ & $70.4$ \\ 
		12 & $\mathbf{403.3}$ & $0.875$ & $654$ & $\mathbf{3.1}$ & $76.0$ &$395.7$ & $0.848$ & $40$ & $5.8$ & $75.6$ &$391.5$ & $0.837$ & $10463$ & $10.7$ & $75.7$ \\ 
		14 & $\mathbf{399.9}$ & $0.860$ & $636$ & $\mathbf{3.5}$ & $80.5$ &$391.2$ & $0.829$ & $41$ & $8.1$ & $79.3$ &$384.5$ & $0.814$ & $14758$ & $22.6$ & $79.4$ \\ 
		16 & $\mathbf{398.7}$ & $0.842$ & $640$ & $\mathbf{3.8}$ & $83.3$ &$388.7$ & $0.808$ & $38$ & $8.2$ & $82.9$ &$385.4$ & $0.800$ & $12004$ & $25.3$ & $82.8$ \\ 
		\hline
	\end{tabular}
	\label{Tab:  SPCA  aver10}
\end{table}
Before reporting the results, we specify the default choices of some parameters in our tests. We choose $\zeta_{\min} = 10^{-20}, \zeta_{\max} = 10^{20}$, $\tau_1=0.999,\tau_2=0.9, \rho=1.5$, and $\lambda={\varepsilon}/{(2R)}$ for RADA-PGD and RADA-RGD. {We also set  $\lambda={\varepsilon}/{(2R)}$ for ARPGDA. } For RADA-RGD, we set $ \nu_{k}=2T_kR^2\beta_{k}, c_1=10^{-4}$, and $ \eta=0.1$. The value of $R$ can be easily computed for each specific problem and is given by $\mu \sqrt{dr}$, $1$, and $\mu N$ for problems \eqref{prob: SPCA ref minmax}, \eqref{prob: FPCA ref minmax}, and \eqref{prob: SSC ref minimax}, respectively. The initial setting of the parameter $\beta_1$ varies across different problem settings. 
For DSGM, ARPGDA, and RADMM, we apply the same line search strategy initialized with the BB stepsize as in RADA-RGD for updating $x$. Unless otherwise stated, we terminate RADA-PGD and RADA-RGD if they return an $\varepsilon$-RGS point. The stopping rules for the other algorithms will be specified later. We utilize the QR factorization and the standard projection as the retraction on the Stiefel manifold and the Grassmann manifold, respectively \cite{absil2008optimization,absil2012projection,bendokat2024grassmann,sato2014optimization}. 
The standard projection onto the Grassmann manifold is obtained by taking the $m$ leading eigenvectors of the matrix to be projected. The Riemannian gradient on the Grassmann manifold is computed as in \cite{sato2014optimization}.
Finally, both the QR-based retraction on the Stiefel manifold and the projection-based retraction on the Grassmann manifold are globally defined for all tangent vectors. Therefore, Assumption \ref{assmption: retraction well defined} holds with $\varrho = +\infty$ for both manifolds.
\subsection{Results on SPCA}

In this subsection, we compare our proposed RADA-RGD with RALM-II (i.e., LS-II in \cite{zhou2023semismooth}) and ManPG-Ada \cite{chen2020proximal}, a more efficient version of ManPG, on synthetic datesets and the real dataset {coil-100} \cite{nene1996columbia}, which contains $7200$ image samples of 100 objects taken from different angles with $d=1024$. For RADA-RGD, we set $\beta_{1}=0.1d\sqrt{r}$ and $T_k\equiv10$. 
The codes for RALM-II and ManPG-Ada are downloaded from \url{https://github.com/miskcoo/almssn} and \url{https://github.com/chenshixiang/ManPG}, respectively. 
For ManPG-Ada, we relax the line search parameter from $0.5$ to $10^{-4}$ and adopt a nonmonotone line search strategy to improve stepsize selection and empirical performance. 
We terminate all three methods once an $\varepsilon$-ROS point is found. 
More specifically, we terminate RADA-RGD when
\begin{equation}\label{equ: epsilon st}
	\max\!\left\{\mathrm{dist}\left(0,\,\mathrm{grad}\,f(x_{k}) + \mathrm{proj}_{\mathrm{T}_{x_k}\mathcal{M}}\left(\nabla \mathcal{A}(x_k)^\top \partial h^*(p_k)\right)\right)\hspace{-2pt},\,\left\|p_k-\mathcal{A}(x_k)\right\|\right\}\leq\varepsilon, 
\end{equation}
where $p_k = y_{k+\frac12}$ is defined in \eqref{y_k+1/2}. For RALM-II, we apply the same criterion \eqref{equ: epsilon st} but with $p_k=y_k$. 
Recall that ManPG solves the subproblem
\begin{align}\label{equ: ManPG subproblem}
	V_k:=\arg\min_{V\in\mathrm{T}_{x_k}\mathcal{M}} \,\left\{\langle \mathrm{grad} f(x_k),V\rangle+\frac{\ell_f}{2}\|V\|^2+h^*(x_k+V) \right\}\nonumber
\end{align}
at each iteration. Here, $\ell_f$ is the Lipschitz constant of $\nabla f$.
We terminate ManPG-Ada if $\max\{\ell_f\|V_k\|,\|V_k\|\}\leq\varepsilon$.
The maximum number of iterations is set as 50000. In each test, we generate the synthetic data matrix $A$ following the method in \cite{zhou2023semismooth} or randomly select 100 samples from the {coil-100} dataset.

First, we set $\varepsilon=10^{-4}$. The average results over 20 runs, each using either a different randomly generated synthetic matrix $A$ or a different set of selected samples from the coil-100 dataset, are reported in Table \ref{Tab:  SPCA aver10}. Here, $\Phi$ is defined in \eqref{Phi}, ``cpu" represents the cpu time in seconds, and “iter" denotes the outer iteration number.
We also compare the normalized variance (denoted by “var") and the sparsity of $X$ (denoted by “spar''). The former is defined as
$\langle AA^\top,\,\hat{X}\hat{X}^\top\rangle/\max_{X\in\mathcal{S}(d,r)}\,\langle AA^\top,\,XX^\top\rangle$, where $\hat{X}$ denotes the returned solution by the algorithm. The latter is defined as the percentage of entries with the absolute value less than $10^{-5}$. 

From Table \ref{Tab:  SPCA  aver10}, we see that the solutions obtained by these algorithms are comparable in terms of the normalized variance and sparsity. Our proposed RADA-RGD always returns the best solution in terms of the value of $\Phi$. 
Moreover, the proposed single-loop RADA-RGD algorithm is significantly more efficient than the compared nested-loop algorithms. 
The performance difference between RADA-RGD and RALM-II is consistent with the analysis in Section \ref{sec: connections}. 
Specifically, although RALM-II requires significantly fewer iterations compared to RADA-RGD, it requires a much higher level of accuracy than RADA-RGD when solving the subproblem. Consequently, RALM-II is less efficient than RADA-RGD for solving the SPCA problems considered.  

\begin{figure}[t]
	\centering
	\begin{subfigure}
		\centering
		\includegraphics[scale=0.4]{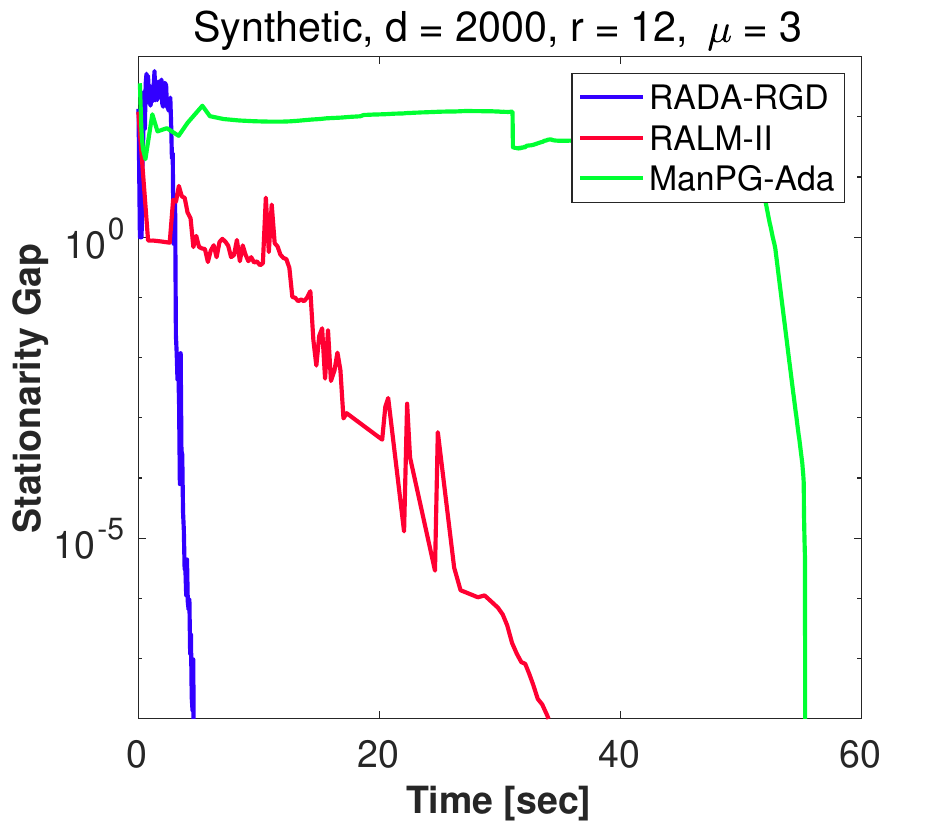}
	\end{subfigure}
	\begin{subfigure}
		\centering
		\includegraphics[scale=0.4]{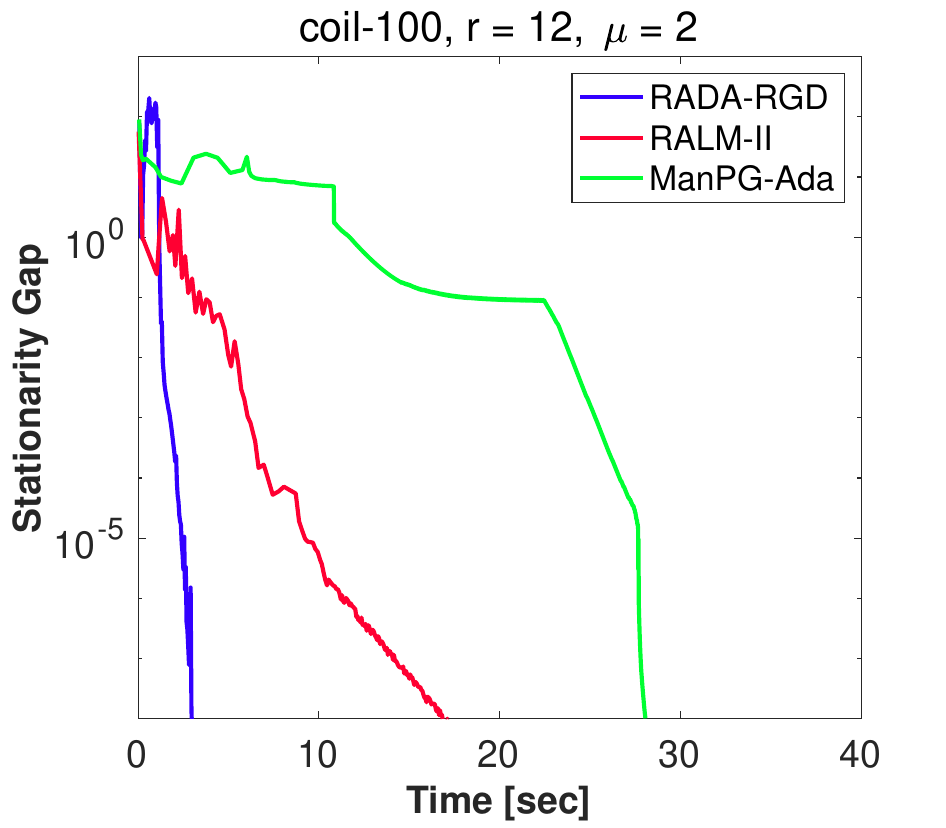}
	\end{subfigure}
	\caption{Stationarity gap versus runtime for two representative instances of SPCA.}
	\label{Fig: Spca}
\end{figure}

Next, we set $\varepsilon = 10^{-8}$ and tighten the subproblem stopping criterion of ManPG-Ada to $10^{-10}$ (otherwise, the stationarity gap stagnates and fails to reach the same accuracy level as RADA-RGD and RALM-II) to observe the stationarity gap versus the runtime of two representative instances in Figure \ref{Fig: Spca}. To ensure a fair and readable comparison, we record the runtime of RALM-II at each outer iteration and that of RADA-RGD and ManPG-Ada every 10 iterations. This is because RALM-II involves the fewest but most expensive outer iterations. As shown in Figure \ref{Fig: Spca}, RADA-RGD reduces the stationarity gap significantly faster than both RALM-II and ManPG-Ada, which is consistent with the results  reported in Table \ref{Tab:  SPCA  aver10}.

\subsection{Results on FPCA}
For the FPCA problem \eqref{prob: FPCA ref minmax},  we generate the data points $\{a_i\}_{i=1}^N$ independently according to the standard Gaussian distribution as in \cite{xu2023efficient2} and divide them into $m=N$ groups. We compare RADA-RGD with the single-loop algorithms ARPGDA \cite{xu2023efficient2}, RADMM \cite{li2023riemannian}, and DSGM \cite{beck2023dynamic}. 
Note that the existing convergence guarantees for RADMM and DSGM do not apply to FPCA since the mapping $\mathcal{A}$ is nonlinear. For RADA-RGD, we set $ \beta_{1}=10^4m^2\sqrt{r}$, $T_k\equiv5$, and $\varepsilon=10^{-8}$.  For ARPGDA, we set $\beta_{k}={10^3m^2\sqrt{r}}/{k^{1.5}}$. For DSGM and RADMM, we set the smoothing parameter $\lambda_k={10}/{k^{1/3}}$. For RADMM, a significant effort is devoted to fine-tuning the penalty parameter in the augmented Lagrangian function, and the choice $\sigma_{k}=10^{-7}k^{1.5}$ is found to work well in our tests. It should be noted that although this choice of $\sigma_k$ violates the convergence condition in \cite{li2023riemannian}, the resulting RADMM demonstrates superior numerical performance.
For the compared algorithms, a more relaxed criterion is employed: The algorithm is terminated if the improvement in $\Phi$ is less than $10^{-8}$ over $1000$ consecutive iterations. The maximum number of iterations for all algorithms is set as 20000.

\begin{figure}[t]
	\centering
	\begin{subfigure}
		\centering
		\includegraphics[scale=0.4]{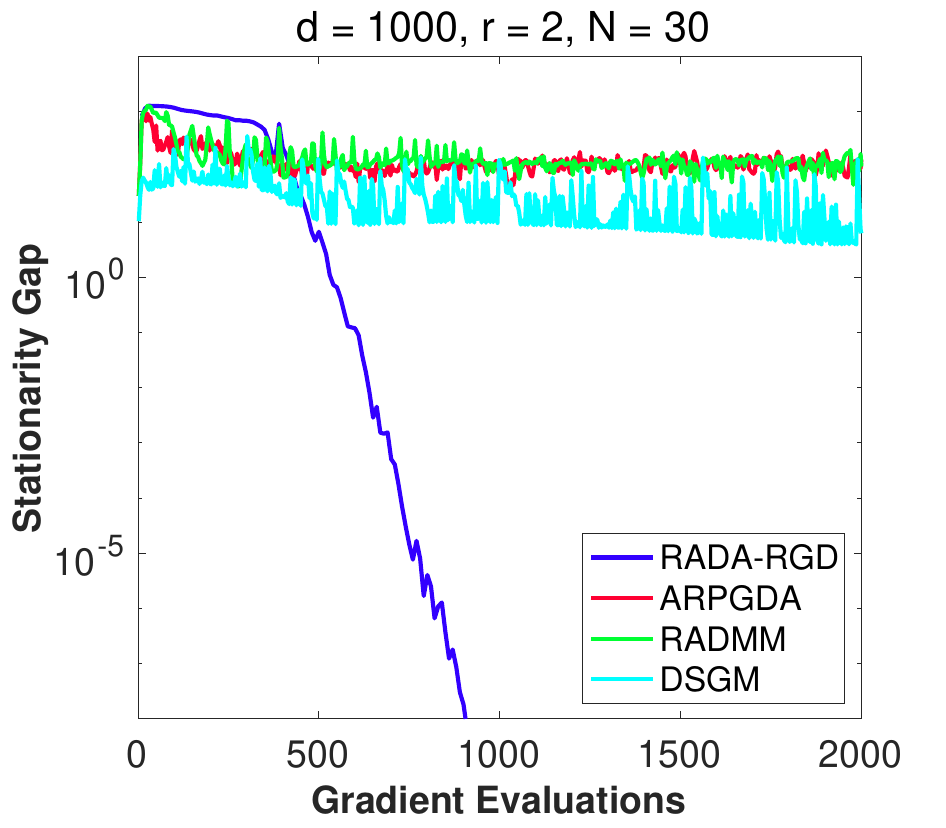}
	\end{subfigure}
	\begin{subfigure}
		\centering
		\includegraphics[scale=0.4]{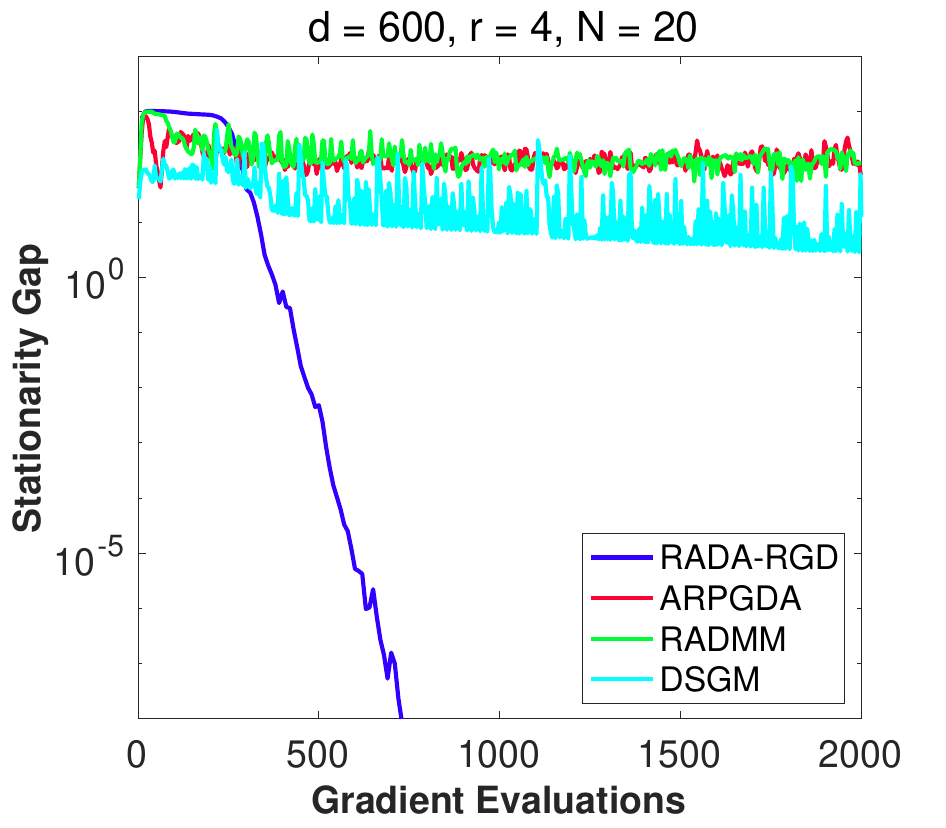}
	\end{subfigure}
	\caption{Stationarity gap versus number of gradient evaluations for two representative instances of FPCA.}
	\label{Fig: fpca}
\end{figure}

\begin{table}
	\centering
	\fontsize{10pt}{\baselineskip}\selectfont
	\caption{{Average performance comparison on FPCA.}}
	\vspace{-0.0cm}
	\tabcolsep=0.15cm
	\renewcommand\arraystretch{1}
	\begin{tabular}{cccc|ccc|ccc|ccc}
		\hline
		& \multicolumn{3}{c}{RADA-RGD} & \multicolumn{3}{c}{ARPGDA} & \multicolumn{3}{c}{RADMM} & \multicolumn{3}{c}{DSGM} \\
		\cline{2-13}
		& $-\Phi$ & cpu & iter & $-\Phi$ & cpu & iter & $-\Phi$ & cpu & iter & $-\Phi$ & cpu & iter \\
		\hline
		$N$ & \multicolumn{12}{c}{$d=1000, r=2$} \\
		\hline
		20 & $\mathbf{120.54}$ & $\mathbf{0.01}$ & $64$ & $118.08$ & $0.24$ & $7715$& $120.18$ & $0.20$ & $6144$& $113.38$ & $0.16$ & $3123$\\
		25 & $\mathbf{98.88}$ & $\mathbf{0.01}$ & $77$ & $97.00$ & $0.24$ & $6708$& $97.94$ & $0.24$ & $6885$& $97.01$ & $0.19$ & $3599$\\
		30 & $\mathbf{84.16}$ & $\mathbf{0.02}$ & $85$ & $82.29$ & $0.22$ & $5764$& $83.11$ & $0.30$ & $7858$& $82.32$ & $0.24$ & $4276$\\
		35 & $\mathbf{73.69}$ & $\mathbf{0.02}$ & $102$ & $72.16$ & $0.33$ & $7833$& $72.80$ & $0.27$ & $6289$& $69.74$ & $0.24$ & $3836$\\
		40 & $\mathbf{65.40}$ & $\mathbf{0.02}$ & $98$ & $64.17$ & $0.32$ & $7257$& $64.77$ & $0.25$ & $5814$& $60.24$ & $0.25$ & $3824$\\
		\hline
		$d$ & \multicolumn{12}{c}{$N=20, r=4$} \\
		\hline
		200 & $\mathbf{56.86}$ & $\mathbf{0.01}$ & $87$ & $56.51$ & $0.02$ & $1579$& $56.59$ & $0.02$ & $1178$& $53.73$ & $0.06$ & $1890$\\
		400 & $\mathbf{104.24}$ & $\mathbf{0.01}$ & $73$ & $103.29$ & $0.07$ & $3174$& $103.08$ & $0.11$ & $4771$& $103.39$ & $0.11$ & $2420$\\
		600 & $\mathbf{146.57}$ & $\mathbf{0.01}$ & $62$ & $145.55$ & $0.32$ & $9975$& $146.20$ & $0.22$ & $6717$& $140.12$ & $0.20$ & $3821$\\
		800& $\mathbf{193.26}$ & $\mathbf{0.01}$ & $53$ & $191.28$ & $0.38$ & $9094$& $192.64$ & $0.31$ & $7248$& $192.25$ & $0.32$ & $5244$\\
		1000 & $\mathbf{237.18}$ & $\mathbf{0.02}$ & $75$ & $235.59$ & $0.51$ & $10124$& $236.65$ & $0.38$ & $7423$& $224.79$ & $0.38$ & $5270$\\
		\hline
	\end{tabular}
	\label{Tab:  FPCA  aver10}
	\label{T}
\end{table}
Table \ref{Tab:  FPCA  aver10} presents the results on FPCA averaged over 20 runs with different randomly generated initial points. From this table, we observe that our proposed RADA-RGD consistently returns the best solutions in terms of the value of $\Phi$ while requiring the least computational time. Specifically, on average, our proposed RADA-RGD is more than 10 times faster than the other compared algorithms when $d$ is large. We further plot the stationarity gap versus the number of gradient evaluations for two representative instances in Figure \ref{Fig: fpca}. As illustrated in the figure, RADA-RGD achieves a faster reduction in the stationarity gap than the other first-order methods. Moreover, the similarities and differences between RADA-RGD and RADMM, as discussed in Section \ref{subsection: Interpreting RADA-RGD from the Perspective of ADMM}, provide important insights into the superior performance of RADA-RGD. In particular, as shown in Algorithm \ref{Algorithm multiBB2}, an additional update of $p$ is performed before updaing $x$, which means that the Riemannian gradient used to update $x$ in RADA-RGD involves the most up-to-date information of $p$. The superior performance of RADA-RGD over ARPGDA further demonstrates that
performing multiple Riemannian gradient steps to ensure sufficient descent of the proposed value function $\Phi_{k}$ is much more efficient than simply performing a Riemannian gradient descent step on certain surrogate of $F(\cdot,y)$ for some fixed $y$, as is done in ARPGDA and its Euclidean counterparts \cite{xu2023unified,pan2021efficient,lu2020hybrid,he2024approximation}.

\subsection{Results on SSC}
In this subsection, we compare our proposed RADA-PGD {and RADA-RGD} with ARPGDA \cite{xu2023efficient2}, DSGM \cite{beck2023dynamic}, and NADMM \cite{lu2018nonconvex} for solving the SSC problem \eqref{probSSC}.  Note that NADMM is specifically designed for solving \eqref{probSSC}.  The main difference between NADMM and RADMM is that NADMM performs a Euclidean projected gradient descent step to update $x$.

First, we perform tests on synthetic datasets, wherein the data points $\{a_i\}_{i=1}^N$ are independently generated according to the standard Gaussian distribution and $W_{ij}= |\langle a_i, a_j\rangle|$ as in \cite{lu2016convex}.  We take $Q_1=X_1X_1^\top$ as the initial point, where $X_1$ consists of the $m$ eigenvectors associated to the $m$ smallest eigenvalues of $L$.  For RADA, we set $\varepsilon=10^{-4}$. For RADA-PGD, we set $\beta_1=N^2\sqrt{m}$, $T_k\equiv1$, and $\zeta_{k,t}={\ell_k^{-1}}$ with $\ell_k={(\lambda+\beta_{k})^{-1}}$. For RADA-RGD, we set $\beta_1=N^2\sqrt{m}$ and $T_k\equiv3$.  For DSGM and NADMM, we set $\lambda_k={10^{-2}}/{k^{1/3}}$.  Additionally, we carefully tune the parameter $\beta_{k}$ in ARPGDA for each test to achieve its best possible performance.  For NADMM, we set the penalty parameter $\sigma_{k}$ in the augmented Lagrangian function as $0.1k^{1/3}$, which enhances the efficiency of NADMM in our tests but does not satisfy the conditions of the theoretical guarantee in \cite{lu2018nonconvex}.  For the compared algorithms,  a more relaxed criterion is employed: The algorithm is terminated if the improvement of the objective function $\Phi$ is less than $10^{-8}$ over $50$ consecutive iterations.

\begin{table}
	\centering
	\fontsize{10pt}{\baselineskip}\selectfont
	\caption{Average performance comparison on SSC (synthetic datasets). 
	}
	\tabcolsep=0.15cm
	\renewcommand\arraystretch{1}
	\begin{tabular}{cccc|ccc|ccc|ccc|ccc}
		\hline
		& \multicolumn{3}{c}{RADA-PGD} &  \multicolumn{3}{c}{RADA-RGD} & \multicolumn{3}{c}{ARPGDA} & \multicolumn{3}{c}{NADMM} & \multicolumn{3}{c}{DSGM} \\
		\cline{2-16}
		& $\Phi$ & cpu & iter & $\Phi$ & cpu & iter & $\Phi$ & cpu & iter & $\Phi$ & cpu & iter & $\Phi$ & cpu & iter \\
		\hline
		$m$ & \multicolumn{15}{c}{$N=200, \mu=0.005$} \\
		\hline
		2 & $\mathbf{1.976}$ & $\mathbf{0.1}$ & $72$ & $\mathbf{1.976}$ & $0.2$ & $67$ & $1.980$ & $0.6$ & $580$& $1.978$ & $0.6$ & $838$& $1.983$ & $8.0$ & $4455$\\
		3 & $\mathbf{2.973}$ & $\mathbf{0.1}$ & $92$ & $\mathbf{2.973}$ & $0.3$ & $78$ & $2.987$ & $6.0$ & $6494$& $2.977$ & $0.8$ & $1136$& $2.988$ & $2.9$ & $2168$\\
		4 & $\mathbf{3.970}$ & $\mathbf{0.2}$ & $104$ & $\mathbf{3.970}$ & $\mathbf{0.2}$ & $77$ & $3.988$ & $5.6$ & $6238$& $3.975$ & $0.9$ & $1324$& $3.989$ & $4.5$ & $3528$\\
		5 & $\mathbf{4.967}$ & $\mathbf{0.2}$ & $123$ & $\mathbf{4.967}$ & $\mathbf{0.2}$ & $76$ & $4.987$ & $6.0$ & $6446$& $4.973$ & $1.2$ & $1759$& $4.992$ & $6.8$ & $5375$\\
		6 & $\mathbf{5.965}$ & $\mathbf{0.2}$ & $141$ & $\mathbf{5.965}$ & $0.3$ & $86$ & $5.987$ & $5.8$ & $6185$& $5.971$ & $1.7$ & $2421$& $5.994$ & $9.8$ & $7601$\\
		\hline
		$\mu$ & \multicolumn{15}{c}{$N=500, m=3$} \\
		\hline
		0.001 & $\mathbf{2.496}$ & $\mathbf{1.0}$ & $116$ & $\mathbf{2.496}$ & $3.0$ & $158$ & $2.500$ & $22.2$ & $3757$& $2.498$ & $3.8$ & $1244$& $2.504$ & $73.5$ & $9198$\\
		0.002 & $\mathbf{2.985}$ & $\mathbf{1.1}$ & $110$ & $\mathbf{2.985}$ & $2.0$ & $93$ & $2.990$ & $46.5$ & $7294$& $2.986$ & $1.9$ & $531$& $2.996$ & $25.9$ & $3084$\\
		0.005 & $\mathbf{3.007}$ & $\mathbf{1.7}$ & $152$ & $\mathbf{3.007}$ & $3.2$ & $143$ & $3.010$ & $3.3$ & $503$& $\mathbf{3.007}$ & $3.0$ & $798$& $3.009$ & $29.7$ & $3282$\\
		0.01 & $\mathbf{3.022}$ & $\mathbf{1.9}$ & $165$ & $\mathbf{3.022}$ & $3.5$ & $152$ & $3.031$ & $12.2$ & $1796$& $\mathbf{3.022}$ & $13.9$ & $3480$& $3.024$ & $16.5$ & $1681$\\
		0.02 & $\mathbf{3.052}$ & $\mathbf{2.6}$ & $227$ & $\mathbf{3.052}$ & $4.1$ & $164$ & $3.057$ & $3.2$ & $436$& $\mathbf{3.052}$ & $17.1$ & $3704$& $\mathbf{3.052}$ & $13.3$ & $1134$\\
		\hline
	\end{tabular}
	\label{Tab:  SSC on Stiefel aver10}
\end{table}
\begin{figure}[t]
	\centering
	\begin{subfigure}
		\centering
		\includegraphics[scale=0.4]{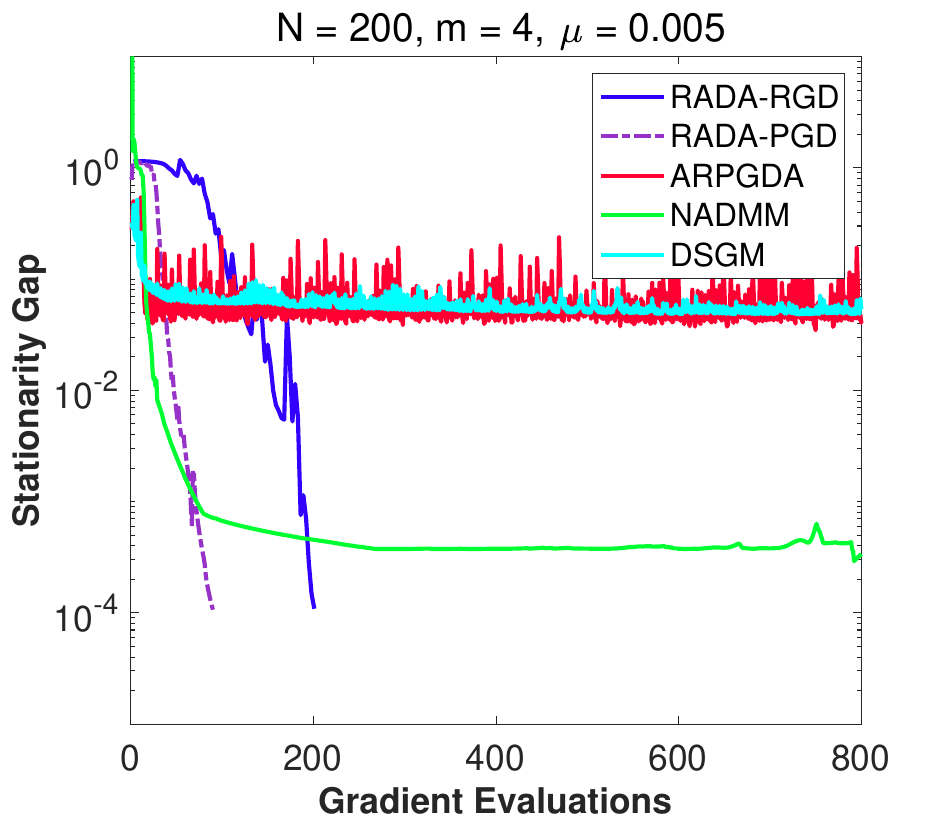}
	\end{subfigure}
	\begin{subfigure}
		\centering
		\includegraphics[scale=0.4]{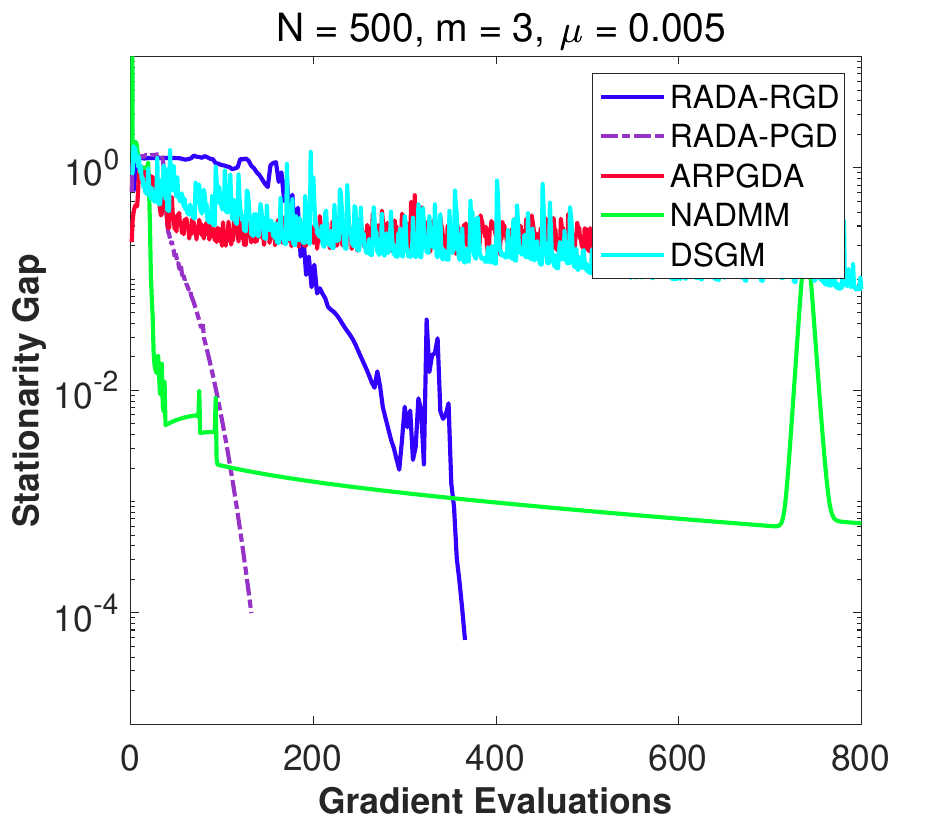}
	\end{subfigure}
	\caption{Stationarity gap versus number of gradient evaluations for two representative instances of SSC.}
	\label{Fig: ssc}
\end{figure}
Table \ref{Tab:  SSC on Stiefel aver10} reports the results on SSC averaged over 20 runs with randomly generated data points, and Figure \ref{Fig: ssc} illustrates the stationarity gap versus the number of gradient evaluations for two representative instances. From both the table and the figure, we observe that the proposed RADA-PGD and RADA-RGD consistently achieve the best values of~$\Phi$ and effectively reduce the stationarity gap.
Notably, RADA-PGD is the most efficient among all compared algorithms. 
Similar to RADA-RGD, RADA-PGD with $T_k\equiv1$ can also be interpreted as a Riemannian sGS-type ADMM. 
This interpretation sheds light on the superior performance of
RADA-PGD over NADMM, given that
sGS-type ADMM has demonstrated better performance than ADMM in many scenarios \cite{chen2017efficient,li2016schur,li2019block}.
Furthermore, RADA-PGD demonstrates even better performance than RADA-RGD, especially when $N=500$. This can be attributed to two reasons. 
First, RADA-PGD only needs to compute the Euclidean gradient instead of the Riemannian gradient, which requires an additional projection onto the tangent space. 
Second, the parameter
$\ell_k$ in \eqref{lemma l-smooth E ineq E} becomes ${(\lambda + \beta_k)^{-1}}$ for this specific problem, making it straightforward to compute. Consequently, we can use ${\ell_k^{-1}}=\lambda+\beta_k$ directly as an appropriate stepsize at the $k$-th iteration without  any additional cost. 

Next, to further evaluate the different algorithms, we compare their clustering performance on four real datasets (i.e., Wine, Iris, Glass, and Letter) from the UCI machine learning repository \cite{uci_repository}. Due to the large size of the Letter dataset, we randomly select a subset of 1000 samples from the first 17  classes. The affinity matrix $W$ is constructed using the Gaussian kernel as in \cite{lu2018nonconvex}, i.e., $ W_{ij}=\mathrm{exp}\left(-{\|a_i-a_j\|^2}/{\kappa}\right)$. We fine-tune the parameters $\kappa$ and $\mu$ for each dataset so that the solutions to the obtained optimization problems yield good clustering performance. We utilize the normalized mutual information (NMI) scores \cite{strehl2002cluster} to measure the clustering performance. Note that a higher NMI score indicates better clustering performance. The parameters and stopping criteria for all algorithms remain the same as those used in the tests on synthetic datasets, with the exception that $\varepsilon=10^{-3}$ is used for both RADA-PGD and RADA-RGD. We also fine-tune the parameter $\beta_{k}$ for ARPGDA in each test. 
\begin{table}[t]
	\centering
	\fontsize{8pt}{\baselineskip}\selectfont
	\caption{Performance comparison on SSC (UCI datasets).} 
	\tabcolsep=0.3cm
	\renewcommand\arraystretch{0.8}
	\begin{tabular}{cccccc}
		\hline
		& RADA-PGD &RADA-RGD & ARPGDA  & DSGM & NADMM \\
		\hline
		\multicolumn{6}{c}{Wine: $N=178, m=3, \kappa=1, \mu=0.001$}\\
		\hline
		NMI & $\mathbf{0.893}$ & $\mathbf{0.893}$ & $\mathbf{0.893}$ & $\mathbf{0.893}$ & $\mathbf{0.893}$ \\ 
		$\Phi$ & $\mathbf{1.899}$ & $\mathbf{1.899}$ & $\mathbf{1.899}$ & $\mathbf{1.899}$ & $\mathbf{1.899}$ \\ 
		iter & $148$ & $115$ & $932$ & $1952$ & $659$ \\ 
		cpu & $\mathbf{0.26}$ & $0.52$ & $1.60$ & $5.12$ & $0.97$ \\ 
		\hline
		\multicolumn{6}{c}{Iris: $N=149, m=3, \kappa=0.2, \mu=0.005$}\\
		\hline
		NMI & $\mathbf{0.776}$ & $\mathbf{0.776}$ & $\mathbf{0.776}$ & $0.757$ & $\mathbf{0.776}$ \\ 
		$\Phi$ & $\mathbf{1.106}$ & $\mathbf{1.106}$ & $\mathbf{1.106}$ & $\mathbf{1.106}$ & $\mathbf{1.106}$ \\ 
		iter & $161$ & $115$ & $320$ & $1491$ & $5995$ \\ 
		cpu & $\mathbf{0.21}$ & $0.44$ & $0.59$ & $3.17$ & $4.82$ \\ 
		\hline
		\multicolumn{6}{c}{Glass: $N=213, m=3, \kappa=1, \mu=0.0002$}\\
		\hline
		NMI & $\mathbf{0.471}$ & $\mathbf{0.471}$ & $\mathbf{0.471}$ & $\mathbf{0.471}$ & $\mathbf{0.471}$ \\ 
		$\Phi$ & $\mathbf{1.818}$ & $\mathbf{1.818}$ & $\mathbf{1.818}$ & $\mathbf{1.818}$ & $\mathbf{1.818}$ \\ 
		iter & $120$ & $63$ & $268$ & $8474$ & $177$ \\ 
		cpu & $\mathbf{0.17}$ & $0.30$ & $0.44$ & $28.35$ & $0.34$ \\ 
		\hline
		\multicolumn{6}{c}{Letter: $N=700, m=12, \kappa=1, \mu=0.001$}\\
		\hline
		NMI & $\mathbf{0.386}$ & $0.384$ & $0.339$ & $0.297$ & $0.315$ \\ 
		$\Phi$ & $\mathbf{16.63}$ & $16.64$ & $16.70$ & $16.73$ & $16.66$ \\ 
		iter & $271$ & $222$ & $4456$ & $10000$ & $10000$ \\ 
		cpu & $\mathbf{8.71}$ & $23.89$ & $216.65$ & $450.76$ & $313.04$ \\ 
		\hline
	\end{tabular}
	\vspace{-0.5cm}
	\label{Tab:  SSC real}
\end{table}

Table \ref{Tab:  SSC real} reports the comparison of NMI scores, objective function values, and computational time for the four real datasets. From the table, we see that RADA-PGD and RADA-RGD consistently find solutions with the best quality in terms of NMI scores and objective function values for all datasets.  For the Letter dataset, the proposed RADA-PGD enjoys a significantly lower computational time, while DSGM and NDAMM hit the preset maximum number of iterations of 10000.

\section{Concluding Remarks}\label{sec: concluding remarks}
In this paper, we introduced a flexible RADA algorithmic framework based on a novel value function for solving a class of Riemannian NC-L minimax problems. Within this framework, we proposed two customized efficient algorithms called RADA-PGD and RADA-RGD, which perform one or multiple projected/Riemannian gradient descent steps and then a proximal gradient ascent step at each iteration. Theoretically, we proved that the proposed algorithmic framework can find both $\varepsilon$-RGS and $\varepsilon$-ROS points of the Riemannian NC-L minimax problem within $\mathcal{O}(\varepsilon^{-3})$ iterations, achieving the best-known iteration complexity. We also showed how the proposed algorithmic framework can be adapted to solve an equivalent Riemannian nonsmooth reformulation of the minimax problem. As a by-product, we made the intriguing observation that RADA and RADA-RGD are closely related to, but essentially different from, RALM and RADMM, respectively. 
Computationally, we demonstrated the superiority of our proposed algorithms compared to existing state-of-the-art algorithms through extensive  numerical results on SPCA, FPCA, and SSC. It is worth noting that our proposed RADA algorithmic framework and the associated algorithms are not derived from existing Euclidean methods; rather, the ideas behind their design even lead to new algorithms for NC-L minimax problems in the Euclidean setting.

Our work suggests several directions for future research. One direction is to design efficient algorithms with strong theoretical guarantees for solving more general Riemannian minimax problems. Another direction is to develop stochastic versions of our algorithms. It is also worthwhile to develop efficient algorithms for broader classes of Riemannian nonsmooth composite problems by suitably exploiting their intrinsic minimax structure.  


%
%
%

\ACKNOWLEDGMENT
{ }
The authors are grateful to the Area Editor, the Associate Editor, and the anonymous reviewers for their valuable comments and suggestions, which have helped improve the quality of the manuscript.
%



\bibliographystyle{informs2014} 
\bibliography{reference_arv} 




%
%
%
%

\end{document}